\documentclass[12pt,dvips]{amsart}
\usepackage{euler, amsfonts, amssymb, latexsym, epsfig,epic}

\newcommand\iso{\cong}
\newtheorem{Theorem}{Theorem}[section]
\newtheorem{Proposition}[Theorem]{Proposition} 
\newtheorem{Lemma}[Theorem]{Lemma}
\newtheorem{Corollary}[Theorem]{Corollary}
\newtheorem*{Corollary*}{Corollary}
\newtheorem*{Theorem*}{Theorem}
\theoremstyle{remark}
\newtheorem*{Definition*}{Definition} 
\newtheorem*{Remark*}{Remark}
\newtheorem{Example}[Theorem]{Example}

        {\begin{list}
                {\noindent\makebox[0mm][r]{\arabic{enumi}.}}
                {\leftmargin=5.5ex \usecounter{enumi}}
        }
        {\end{list}}

\newenvironment{romanlist}%
        {\begin{list}
                {\noindent\makebox[0mm][r]{(\rm\roman{enumi})}}
                {\leftmargin=5.5ex \usecounter{enumi}}
        }
        {\end{list}}

\newenvironment{alphlist}%
        {\begin{list}
                {\noindent\makebox[0mm][r]{\rm\alph{enumi})}}
                {\leftmargin=5.5ex \usecounter{enumi}}
        }
        {\end{list}}

\newcommand\integers{{\mathbb Z}}

\newcommand\naturals{{\mathbb N}}
\newcommand\Spec{{\rm Spec\:}}

\theoremstyle{plain}

\newcommand\IN{\mathsf{in}}
\newcommand\init{\mathsf{in}}
\newcommand\val{\mathsf{val}}

\newcommand\Le{\prec}
\newcommand\FT{\mathit{FT}}
\newcommand\FST{\mathit{FSVT}}
\newcommand\SVT{\mathit{SVT}}
\newcommand\LCM{\mathrm{LCM}}
\newcommand\Ess{{\mathcal E\hspace{-.2ex}s\hspace{-.1ex}s}}
\newcommand\rr[3]{\makebox[0pt][l]{$\phantom{r}^{#3}$}r_{#1#2}}
\newcommand\rrr[1]{r^{#1}}
\newcommand\tpi{{\widetilde\pi}}

\def\dfn{\textbf} % maybe should be {\emph}?
\def\sub#1#2{_{#1 \times #2}}
\def\<{\langle}
\def\>{\rangle}
\def\minus{\smallsetminus}

\newcommand{\excise}[1]{}%{$\star$\textsc{#1}$\star$}
\newcommand{\comment}[1]{$\star${\sf\textbf{#1}}$\star$}
\def\GG{{\mathcal G}}
\def\PP{{\mathcal P}}
\def\SS{{\mathfrak S}}
\def\ZZ{{\mathbb Z}}

\def\kk{\Bbbk}
\def\xx{\mathbf{x}}
\def\yy{\mathbf{y}}
\def\zz{\mathbf{z}}
\def\wt{\mathit{wt}}

\def\ol#1{{\overline {#1}}}

\newenvironment{proofof}[1]{\begin{trivlist}\item {\it
        Proof of {#1}.\,}}{\mbox{}~\hfill$\Box$\end{trivlist}%
        \mbox{}\vspace{-1ex}}

\font\co=lcircle10
\def\petit#1{{\scriptstyle #1}}

\def\jr{\smash{\raise2pt\hbox{\co \rlap{\rlap{\char'005} \char'007}}
               \raise6pt\hbox{\rlap{\vrule height5pt}}
               \raise2pt\hbox{\rlap{\hskip4pt \vrule height0.4pt depth0pt
                width5.7pt}}
               \raise2pt\hbox{\rlap{\hskip-9.5pt \vrule height.4pt depth0pt
                width6.2pt}}
               \lower6pt\hbox{\rlap{\vrule height4.5pt}}}}
\def\rj{\smash{\raise2pt\hbox{\co \rlap{\rlap{\char'004} \char'006}}
               \raise6pt\hbox{\rlap{\vrule height5pt}}
               \raise2pt\hbox{\rlap{\hskip4pt \vrule height0.4pt depth0pt
                width5.7pt}}
               \raise2pt\hbox{\rlap{\hskip-9.5pt \vrule height.4pt depth0pt
                width6.2pt}}
               \lower6pt\hbox{\rlap{\vrule height4.5pt}}}}
\def\je{\smash{\raise2pt\hbox{\co \rlap{\rlap{\char'005}
                \phantom{\char'007}}}\raise6pt\hbox{\rlap{\vrule height5pt}}
               \raise2pt\hbox{\rlap{\hskip-9.5pt \vrule height.4pt depth0pt
                width6.2pt}}}}
\def\ej{\smash{\raise2pt\hbox{\co \rlap{\rlap{\char'004}\phantom{\char'006}}}
               \raise2pt\hbox{\rlap{\hskip-9.5pt \vrule height.4pt depth0pt
                width6.2pt}}
               \lower6pt\hbox{\rlap{\vrule height4.5pt}}}}
\def\er{\smash{\raise2pt\hbox{\co \rlap{\rlap{\phantom{\char'005}} \char'007}}
               \raise2pt\hbox{\rlap{\hskip4pt \vrule height0.4pt depth0pt
                width5.7pt}}
               \lower6pt\hbox{\rlap{\vrule height4.5pt}}}}
\def\re{\smash{\raise2pt\hbox{\co \rlap{\rlap{\phantom{\char'004}} \char'006}}
               \raise6pt\hbox{\rlap{\vrule height5pt}}
               \raise2pt\hbox{\rlap{\hskip4pt \vrule height0.4pt depth0pt
                width5.7pt}}}}
\def\+{\smash{\lower6pt\hbox{\rlap{\vrule height17pt}}
                \raise2pt%%%% for ordinary pipe dreams
                \hbox{\rlap{\hskip-9pt \vrule height.4pt depth0pt
                width18.7pt}}}}
\def\hor{\smash{\raise2pt\hbox{\rlap{\hskip-9.5pt \vrule height.4pt depth0pt
                width19.2pt}}}}
\def\ver{\smash{\lower6pt\hbox{\rlap{\vrule height17pt}}}}
\def\ho{\smash{\hbox{\rlap{\vrule height5pt}}
                \raise2pt%%%% for ordinary pipe dreams
                \hbox{\rlap{\hskip-9pt \vrule height.4pt depth0pt
                width18.7pt}}}}

\def\perm#1#2{\hbox{\rlap{$\petit {#1}_{\scriptscriptstyle #2}$}}%
                \phantom{\petit 1}}

\def\textcross{\ \smash{\lower4pt\hbox{\rlap{\hskip4.15pt\vrule height14pt}}
                \raise2.8pt\hbox{\rlap{\hskip-3pt \vrule height.4pt depth0pt
                width14.7pt}}}\hskip12.7pt}
\def\textelbow{\ \hskip.1pt\smash{\raise2.75pt%
                \hbox{\co \hskip 4.15pt\rlap{\rlap{\char'004} \char'006}
                \lower6.8pt\rlap{\vrule height3.5pt}
                \raise3.6pt\rlap{\vrule height3.5pt}}
                \raise2.8pt\hbox{%
                  \rlap{\hskip-7.15pt \vrule height.4pt depth0pt width3.5pt}%
                  \rlap{\hskip4.05pt \vrule height.4pt depth0pt width3.5pt}}}
                \hskip8.7pt}

\newcommand{\cellsize}{20}
\newlength{\cellsz} 
\newsavebox{\cell}
\sbox{\cell}{\begin{picture}(\cellsize,\cellsize)
\put(0,0){\line(1,0){\cellsize}}
\put(0,0){\line(0,1){\cellsize}}
\put(\cellsize,0){\line(0,1){\cellsize}}
\put(0,\cellsize){\line(1,0){\cellsize}}
\end{picture}}
\newcommand\cellify[1]{\def\thearg{#1}\def\nothing{}%
\ifx\thearg\nothing
\vrule width0pt height\cellsz depth0pt\else
\hbox to 0pt{\usebox{\cell} \hss}\fi%
\vbox to \cellsz{
\vss
\hbox to \cellsz{\hss$#1$\hss}
\vss}}
\newcommand\tableau[1]{\vtop{\let\\\cr
\baselineskip -16000pt \lineskiplimit 16000pt \lineskip 0pt
\ialign{&\cellify{##}\cr#1\crcr}}}

\begin{document}%%%%%%%%%%%%%%%%%%%%%%%%%%%%%%%%%%%%%%%%%%%%%%%%%%%%%%
%%%%%%%%%%%%%%%%%%%%%%%%%%%%%%%%%%%%%%%%%%%%%%%%%%%%%%%%%%%%%%%%%%%%%%

\pagestyle{plain}

\mbox{}\vspace{-2.6ex}
\title{Gr\"obner geometry of vertex decompositions \\
  and of flagged tableaux\vspace{-1ex}}
\author{Allen Knutson}
%\thanks{AK was supported by an NSF grant.}
\address{Department of Mathematics, University of California,
        Berkeley, CA 94720, USA}
\email{allenk@math.berkeley.edu}
\author{Ezra Miller\vspace{-1ex}}
%\thanks{AK was supported by an NSF grant and
%EM was supported by NSF grant DMS-0304789.}
\address{School of Mathematics, University of Minnesota,
        Minneapolis, MN 55455, USA} 
\email{ezra@math.umn.edu}
\author{Alexander Yong}
\address{Department of Mathematics, University of California,
        Berkeley, CA 94720, USA
        \indent{\itshape and}
        The Fields Institute, 222 College Street, Toronto, Ontario,
        M5T 3J1, Canada}
% \begingroup{\nobreak\indent{\itshape and}{,
% \ignorespaces##1\unskip}\/:\space   ##2\par}\endgroup
\email{ayong@math.berkeley.edu}
%\date{\today}
\date{7 March 2005}

\begin{abstract}
  We relate a classic algebro-geometric degeneration technique, dating
  at least to [Hodge 1941], to the notion of vertex decompositions of
  simplicial complexes.  The good case is when the degeneration is
  reduced, and we call this a \emph{geometric vertex decomposition}.
  
  Our main example in this paper is the family of \emph{vexillary matrix
    Schubert varieties}, whose ideals are also known as (one-sided) ladder
  determinantal ideals.  Using a diagonal term order to specify the
  (Gr\"obner) degeneration, we show that these have geometric vertex
  decompositions into simpler varieties of the same type.  {}From this,
  together with the combinatorics of the pipe dreams of [Fomin--Kirillov
  1996], we derive a new formula for the numerators of their multigraded
  Hilbert series, the double Grothendieck polynomials, in terms of
  \emph{flagged set-valued tableaux}. This unifies work of [Wachs 1985]
  on flagged tableaux, and [Buch 2002] on set-valued tableaux, giving
  geometric meaning to both.
% \comment{EM: Buch will complain that his paper actually does the
% generalization to two sets of variables, as he does mention it in a
% remark.  That's okay, though; it will be clear from the word
% ``double'' before ``Grothendieck polynomials'' that we work in two
% sets of variables.} 
% Moreover, our approach naturally produces the generalization to two
% sets of variables.
  
  This work focuses on diagonal term orders, giving results
  complementary to those of [Knutson--Miller 2005], where it was shown
  that the generating minors form a Gr\"obner basis for any
  \emph{anti}diagonal term order and \emph{any} matrix Schubert variety.
  We show here that under a diagonal term order, the only matrix
  Schubert varieties for which these minors form Gr\"obner bases are the
  vexillary ones, reaching an end toward which the ladder determinantal
  literature had been building.
\end{abstract}

\maketitle

\vspace{-5ex}
{\small \tableofcontents}

%%%%%%%%%%%%%%%%%%%%%%%%%%%%%%%%%%%%%%%%%%%%%%%%%%%%%%%%%%%%%%%%%%%%%%
\section{Introduction and statement of results}%%%%%%%%%%%%%%%%%%%%%%%
%%%%%%%%%%%%%%%%%%%%%%%%%%%%%%%%%%%%%%%%%%%%%%%%%%%%%%%%%%%%%%%%%%%%%%

Fix an ideal~$I$ in a polynomial ring, or correspondingly,
its zero scheme $X$ inside a coordinatized vector space.  
Each term order yields a
Gr\"obner basis for~$I$, or geometrically, a Gr\"obner degeneration of
$X$ into a possibly nonreduced union of coordinate
subspaces.  Such a degeneration often creates too many pieces all at
once, or spoils geometric properties like reducedness; 
it can be better to work instead with less drastic degenerations
that take the limit of~$X$ under rescaling just one axis at a time.
The limit~$X'$ breaks into two collections of pieces: a
\emph{projection part}\/ and a \emph{cone part}.  In cases where
$X'$ is reduced, quantitative information such as
multidegrees and Hilbert series of the original variety~$X$ can be
derived separately from the parts of this \emph{geometric vertex
decomposition}\/ of~$X$ and combined later.  Reducing the computation
of invariants of~$X$ to those of~$X'$ can be especially
helpful when the projection and cone parts of~$X'$ are simpler
than~$X$.

Under suitable hypotheses, repeating the degeneration-decomposition
procedure for each coordinate axis in turn eventually yields the
Gr\"obner degeneration, but with extra inductive information.  When
the limit $X''$ of this sequence is a union of coordinate subspaces,
or equivalently, $X''$ is defined by a squarefree monomial ideal, the
inductive procedure corresponds exactly to the usual notion of vertex
decomposition for simplical complexes, as defined in~\cite{BP}.

Our goals in this paper are to introduce and develop foundations of
% the algebra of 
geometric vertex decompositions, to apply these generalities to the
class of vexillary matrix Schubert varieties, and to exhibit the
resulting combinatorics on their Gr\"obner degenerations for diagonal
term orders.  In particular, through the notion of \emph{flagged
set-valued tableaux}, we unify the work of Wachs on flagged tableaux
\cite{Wachs} and Buch on set-valued tableaux \cite{Buch}, giving
geometric meaning to both. Moreover, using these tableaux, we obtain
new formulae for homological invariants of the vexillary matrix
Schubert varieties.  Our results can be interpreted as providing a
complete, general, combinatorially enriched development of the theory
surrounding Gr\"obner bases for the extensively studied
\emph{ladder determinantal ideals},
which are the defining ideals of vexillary matrix Schubert varieties.

We begin in this section with a more precise overview, including
statements of our main theorems.

\subsection{Vertex decompositions of simplicial complexes}\label{sub:vd}

Let $\Delta$ be a simplicial complex, and $l$ a vertex of~$\Delta$.
Define two subcomplexes of~$\Delta$: the \dfn{deletion} of~$l$ is the
union~$\delta$ of faces of~$\Delta$ not containing~$l$, and the
\dfn{star} of~$l$ is the union~$\sigma$ of the closed faces
containing~$l$.  Then $\Delta$ equals the union of~$\delta$
and~$\sigma$ along the \dfn{link} $\lambda = \delta \cap \sigma$
of~$l$.  The expression $\Delta = \delta \cup_\lambda \sigma$
of~$\Delta$ as a union of the deletion and star of~$l$ glued along the
link is called the \dfn{vertex decomposition} of the complex~$\Delta$
at the vertex~$l$.  Note that $\sigma$ is a cone, namely the cone on
$\lambda$ with cone point~$l$.

Vertex decompositions allow for inductive calculations on simplicial
complexes, deriving good properties of~$\Delta$ from corresponding
properties of~$\delta$ and~$\lambda$.  One such property is
shellability, as first related to vertex decompositions in~\cite{BP}.

\subsection{An analogue for affine schemes}\label{ssec:gvd}

Associated to a simplicial complex $\Delta$ with vertex set~$V$ is the
\dfn{Stanley--Reisner scheme} $\Spec \kk[\Delta]$, a reduced scheme in
the vector space ${\mathbb A}^{\!V}$ over~$\kk$ with basis~$V$ defined
by
\[
  \Spec \kk[\Delta] = \bigcup_{F \in \Delta} {\mathbb A}^{\!F},
% \, \{\hbox{the coordinate subspace of~$V$ with basis }F\}.
\]
where ${\mathbb A}^{\!F}$ is the coordinate subspace of~${\mathbb
A}^{\!V}$ with basis~$F$.  {}Starting from this perspective, the
notion of vertex decomposition extends to any coordinatized affine
scheme~$X$, meaning a subscheme of the vector space $H \times L$
where $L$ is a line with coordinate~$y$.
% In the case of $X = \Spec \kk[\Delta]$, this will reproduce the
% definition in Section~\ref{sub:vd}.  One new complication is that
% $X$ may be irreducible, so our first step is to degenerate it to
% something reducible.
As $X$ might be irreducible, part of the extension involves
breaking~$X$ into pieces by degeneration.

On the vector space $H \times L$ we have an action of the algebraic
torus $\kk^\times = \kk \minus \{0\}$ by scaling the~$L$ coordinate:
$t \cdot (\vec x,y) = (\vec x,t y)$.  Consider the flat limit $X' =
\lim_{t\to 0} t\cdot X$, which is the result of a sort of gradual
projection of~$X$ to~$H$. % = H \times \{0\}$.

That $X'$ contains the closure~$\Pi$ of the actual projection of~$X$
to~$H$ is obvious;
%, which will usually have the same dimension as $X$ itself.
but there is usually more in~$X'$.  View $L$ as the finite part of $L
\cup \{\infty\} \cong {\mathbb P}^1$, so $H \times L \subset H \times
{\mathbb P}^1$, and we may take the closure~$\ol X$ of~$X$ inside $H
\times {\mathbb P}^1$.  As we slow-motion project~$X$ to $H \times
\{0\}$ by scaling the $L$ coordinate, we pull it away from $H_\infty =
H \times \{\infty\}$ without changing the intersection $\Lambda = \ol
X \cap H_\infty$.
% $\subseteq H \times \{\infty\} \iso H$.
Consequently, the limit~$X'$ must contain the cone $\Lambda \times L$
from the origin over~$\Lambda$.
$$X' \supseteq (\Pi\times \{0\}) \cup_{\Lambda\times \{0\}} (\Lambda\times L)$$

By part of our first geometric result, Theorem~\ref{thm:geoGVD}, this
lower bound correctly calculates $X'$ as a \emph{set}.
%\mbox{$\Pi \cup_{\Lambda\!}  (\Lambda\times L)$}.  
When $X'$ and this union are equal \emph{as schemes}\/ and not just as
sets, such as when the ideal $I(X')$ is radical and hence the
intersection $P \cap C$ of radical ideals $P = I(\Pi)$ and $C =
I(\Lambda \times L)$, we call this the \dfn{geometric vertex
decomposition of~$X$ along the splitting $H \oplus L$}.  The letter
$P$ here stands for ``projection'', while $C$ stands for ``cone''.

\begin{Example} \label{ex:hyperbola}
Let $I = \<xy-1\>$, so $X$ is a hyperbola.  Its projection to the
$x$-axis is dense, and its intersection with the line at infinity
meeting the $y$-axis occurs at $\infty$ on the $y$-axis itself.
% , at $x=0$.
Hence we expect the limit~$X'$ to contain the $x$-axis and the
$y$-axis.  Calculating the limit using the $\kk^\times$ action $t\cdot
I = \<xy-t\>$, we set $t = 0$ to get the ideal of~$X'$, namely $I' =
\<xy\>$.  The equality $I' = P \cap C$ holds here, as $P = \<y\>$ and
$C = \<x\>$, so this example is a geometric vertex decomposition of
the hyperbola. The ideal $I'$ is Stanley--Reisner, but fairly 
dull---the simplicial complex consists of two points.
\end{Example}

\begin{Example}\label{ex:hyperbolaplust}
Let $I=\< xy-1 \> \cap \< x,y\>$, so $X$ is now a hyperbola union a
point at the origin.  A Gr\"obner basis for $I$ under either of the
lexicographic term orders is given by $\{y^{2}x-y, yx^{2}-x\}$.  Thus
$I' = \< y^2 x, yx^2\>$.  As $P=\<y\>$ and $C=\<x\>$, we get $I' \neq
C \cap P = \<xy\>$, so we do not have a geometric vertex decomposition
in this case; the limit~$X'$ contains more, scheme-theoretically, than
just the union~of~two~lines. Note that even though $P,C$ are radical ideals,
$I'$ is not radical.
\end{Example}

\begin{Example}\label{ex:SR}
Let $\Delta$ be a simplicial complex and $X=\Spec \kk[\Delta]$ be the 
associated Stanley--Reisner scheme. Let $L$ be the line corresponding
to a vertex $l$ and $H$ be the hyperplane defined by the sum of the other
coordinates. Then the limiting process does not
change $X$, i.e., $X'=X$, and the geometric vertex decomposition
\[
  X' = (\Pi\times \{0\})\cup_{\Lambda\times \{0\}}(\Lambda\times L)
\]
is obtained by applying the Stanley--Reisner recipe 
to each of the subcomplexes of the vertex decomposition
\[
  \Delta = \delta\cup_\lambda \sigma.
\]
\end{Example}

The algebro-geometric degeneration technique investigated here is
classical, see, e.g., \cite{Hodge41}.  However, our desire to
explicate the analogy with vertex decompositions of simplicial
complexes, and to put this story inside a general framework,
was motivated by our work with vexillary matrix Schubert varieties, as
detailed below.  Actually, during the preparation of this text,
further examples of geometric vertex decompositions and
applications of the methods given below have been found, see, 
e.g.,~\cite{KYI, KYII, Ksubword, KMN, KZJ, PS} 
(see also~\cite{K:slides} for a summary of the first four of these papers).  
We believe that it would be very interesting to find more examples of
schemes that can be profitably studied from the viewpoint suggested here.

\subsection{Vexillary matrix Schubert varieties} \label{sub:vex}

To each permutation $\pi\in S_n$, there is an associated 
\emph{matrix Schubert variety}\/ $\ol X_\pi \subseteq M_n$ living in the
space~$M_n$ of $n\times n$ matrices,
% (this definition from \cite{fulton:92} is recapitulated in
% Section~\ref{sec:gvdvmSv}),
and also a set of \emph{accessible boxes}\/ in the \emph{diagram}\/
of~$\pi$ (these definitions appear in
Sections~\ref{sub:matrixSchub}--\ref{ssec:vexillary}).  Each
accessible box $(l,m)$ yields a splitting $M_n = H\oplus L$, where $H$
consists of the matrices with entry~$0$ at $(l,m)$, and $L$ consists
of those matrices with entries~$0$ everywhere else.

We are now ready to attempt a geometric vertex decomposition of $\ol
X_\pi$.  Things behave particularly well if $\pi$ is \emph{vexillary},
a condition with surprisingly many equivalent formulations, some of
which are expounded in Section~\ref{ssec:vexillary}.\footnote{%
  See \cite{KYII} for a treatment of the non-vexillary case.} Among the
vexillary permutations, the simplest (and easiest to characterize) are
the \dfn{Grassmannian permutations}, which are those with exactly one
descent $\pi(i)>\pi(i+1)$.  The following is stated more precisely in
Lemma~\ref{lemma:chain}, Theorem~\ref{thm:grobnervex}, and
Proposition~\ref{p:inductionstep}.

\begin{Theorem*}
Fix a vexillary permutation~$\pi$.
  \begin{alphlist}
  \item 
  Let $M_n = H\oplus L$ be the decomposition at an accessible box.
  Set $X' = \lim_{t\to 0} t\cdot \ol X_\pi$.  Then $X$ has a geometric
  vertex decomposition given by 
  $X' = (\Pi \times \{0\}) \cup_{\Lambda\times \{0\}} (\Lambda\times L)$, 
  where $\Pi \times L$ and $\Lambda\times L$ are each again vexillary
  matrix Schubert varieties.
% (up to crossing with $L$ \comment{EM: have I got it right which one
% requires the crossing with $L$?})
% 
% EM: the following sentence doesn't add much to this theorem (and
% essential minors haven't been introduced before this point):
% Moreover, this limit can be computed from the essential minors for
% $\ol X_\pi$. 
  \item
  There is a sequence of permutations $\sigma_1,\ldots,\sigma_t$ with
  $\sigma_1$ being Grassmannian and \mbox{$\sigma_t = \pi$}, such that
  each $\ol X_{\sigma_i}$ for $i > 1$ arises as the cone part of a
  geometric vertex decomposition of the previous $\ol X_{\sigma_{i-1}}$.
  \end{alphlist}
\end{Theorem*}

One interpretation of these two statements is that, if we wish to
inductively study Grassmannian matrix Schubert varieties via our
degeneration technique, then the class of vexillary matrix Schubert
varieties is a suitable one to work in.  Another interpretation is
that if we wish to answer questions about vexillary matrix Schubert
varieties, then we should try to reduce these to questions about
Grassmannian matrix Schubert varieties.

\subsection{Gr\"obner bases} \label{sub:gb}

% EM: The following is covered before Section 1.1 now:
% The algebraic way to study the end result of this inductive process
% is with Gr\"obner bases.
% %where the decomposition $H\oplus L$ leads to a term order in which
% %the coordinate $y$ on $L$ is expensive. 
% We will be interested in \dfn{diagonal} term orders, in which the
% leading term of any minor is its diagonal term.

The ideal $I_\pi$ of any matrix Schubert variety~$\ol X_\pi$ is
generated by certain minors in the $n \times n$ matrix of variables, 
namely Fulton's \emph{essential minors}\/ \cite{fulton:92}.  In this
paper we are concerned with \dfn{diagonal} term orders, which by
definition choose from each minor its diagonal term as the largest.
The following combines Theorem~\ref{thm:grobnervex} and
Theorem~\ref{thm:neg} into one statement.

\begin{Theorem*}
%   \begin{itemize}
%   \item
%     If $\pi$ is vexillary, then the essential minors are a
%     Gr\"obner basis for $\ol X_\pi$.
%   \item
% EM: The first sentence of the following doesn't imply the second:
%     If $\pi$ is not vexillary, then there are two essential
%     determinants for $\ol X_\pi$ such that the diagonal term of one
%     divides that of the other.  In particular the essential
%     determinants are not a Gr\"obner basis.
%   \end{itemize}
Let $\pi$ be a permutation.  The essential minors constitute a
Gr\"obner basis for $\ol X_\pi$ under some (and hence any) diagonal
term order if and only if $\pi$ is vexillary.
\end{Theorem*}

The commutative algebra literature has been for some time sneaking up
on the ``if'' direction of this result: extensively studied 
classes of increasingly complex \emph{ladder determinantal 
ideals}\/ defined over the past
decade are special cases of vexillary Schubert determinantal ideals;
see, e.g., \cite[Section~2.4]{grobGeom} \cite{Stu90, fulton:92, GM00}
and the references therein.   
The most inclusive class of ladder
determinantal ideals whose generating minors have been shown
previously to form diagonal Gr\"obner bases appear in \cite{GM00}
(which also contains a well-written exposition about past
developments) and cover
% almost (but not quite) all
a substantial portion of the vexillary~cases.

The ``only if'' direction is striking because the essential minors
in~$I_\pi$ form a Gr\"obner basis for any \emph{antidiagonal}\/ term
order, even if the permutation $\pi$ is not vexillary \cite{grobGeom}.
%This means an order that selects the antidiagonal term from the $k!$
%terms of each $k\times k$ minor.
Hence we push the diagonal term orders as far as they can go.  The key
point is that a permutation fails to be vexillary precisely when two of its
essential rank conditions are nested, causing the diagonal terms of
some essential minors to divide the diagonal terms of other (larger)
essential minors.

\subsection{Flagged set-valued tableaux} \label{sub:fsvt}

As promised earlier, the initial scheme that is produced in
Section~\ref{sub:gb} by Gr\"obner-degenerating the matrix Schubert
variety~$\ol X_\pi$ all at once exhibits inductive combinatorial
structures inherited from stepwise geometric vertex decompositions.
%Surprisingly, the flavor---and also many details---of the results
%parallel those of \cite{grobGeom}, even though the term orders there
%were antidiagonal while those considered here are diagonal.

In \cite{grobGeom}, the antidiagonal initial schemes of all matrix Schubert
varieties were shown to be Stanley--Reisner schemes of certain
\emph{subword complexes}\/ (whose definition from \cite{subword}
%in the generality of arbitrary words in~$S_n$ 
we recall in Section~\ref{section:subword}, in a special case).
The faces of the subword complexes in
\cite{grobGeom} corresponded naturally to the \emph{reduced pipe
dreams}\/ of Fomin and Kirillov \cite{FKyangBax,BB}.

Here,
% where the term orders are diagonal and the permutations vexillary,
we
% (perhaps surprisingly)
again get subword complexes for the initial scheme, along with a
geometric explanation for their vertex-decomposability.  However, the
combinatorics involves \emph{flagged set-valued tableaux}, whose
definition we introduce in Section~\ref{sec:fsvt}, providing the
natural common generalization of \emph{flagged tableaux}\/ \cite{Wachs}
and \emph{set-valued tableaux}\/ \cite{Buch}.\footnote{An alternative
combinatorial approach to this simplicial complex structure on
tableaux, logically independent of~\cite{subword}, is given
in~\cite{KMYII}.}  For the term order, pick a total ordering of the
$n^2$ matrix variables in which no variable appears earlier in the
order than another one weakly to the southeast; the resulting
lexicographic term order is easily seen to be a diagonal term
order. The following theorem provides our Gr\"{o}bner geometry
%/Gr\"{o}bner degeneration 
explanation of the naturality of flagged
set-valued tableaux.

\begin{Theorem*}
  If $\pi$ is a vexillary permutation, then the above lex ordering
  induces a sequence of degenerations of~$\ol X_\pi$, each one a
  geometric vertex decomposition.  The end result of these
  degenerations is a vertex decomposition of the initial scheme, the
  Stanley--Reisner scheme of a certain subword complex~$\Gamma_\pi$,
  whose interior faces correspond in a natural way to flagged
  set-valued tableaux and whose maximal faces correspond to flagged
  tableaux.
\end{Theorem*}

The first sentence of the above result was discussed in Section~\ref{sub:vex}.
To relate a geometric vertex decomposition along the way to an actual
vertex decomposition at the end of the sequence we use 
Proposition~\ref{prop:inductivegvd}.
The rest is a combination of
Theorem~\ref{thm:grobnervex'}, where we identify the initial ideal,
and Theorem~\ref{t:FSVT}, where we biject the flagged set-valued
tableaux with special cases of combinatorial diagrams called \emph{pipe
dreams}.  When the permutation $\pi$ is Grassmannian, the interior
faces correspond bijectively to the set-valued tableaux---with no
flagging---for the associated partition (Theorem~\ref{t:SVT}), and the
facets correspond to the usual semi\-standard Young tableaux
(Proposition~\ref{p:SSYT}).  The subword combinatorics of
\mbox{diagonal} initial ideals is new even for the ladder
determinantal ideals whose Gr\"obner bases were already
known~from~\cite{GM00}.

\subsection{Schubert and Grothendieck polynomials} \label{sub:schub}

The combinatorics of initial ideals yields formulae for homological
invariants.
% such as multidegrees and Hilbert series.

Geometrically, the space $M_n$ of $n \times n$ matrices carries two
actions of the group~$T$ of invertible diagonal matrices, by
multiplication on the left and inverse multiplication on the right,
each preserving every decomposition $H\oplus L$ along a matrix entry.
The resulting Gr\"obner degenerations are therefore $T\times
T$-equivariant, so they preserve the
% $K_{T\times T}$-polynomial and $T\times T$-multidegree
% % EM: The appendix mentioned in the next line no longer exists:
% (see the appendix to \cite{grobGeom} for these definitions).
$T\times T$-equivariant classes of~$\ol X_\pi$ in both cohomology and
\mbox{$K$-theory}.

Equivalently (and more algebraically), the coordinate ring of~$\ol
X_\pi$ carries a grading by~$\ZZ^{2n}$
% $ = \bigoplus_{i=1}^n \ZZ \cdot \{x_i,y_i\}$
in which the variable~$z_{ij}$ has ordinary weight $x_i - y_j$.  Under
this grading, the variety $\ol X_\pi$ has the same $\ZZ^{2n}$-graded
$\mathit{K}$-polynomial and multidegree (see \cite[Chapter~8]{cca} for
definitions) as the Stanley--Reisner initial scheme from
% the previous subsection
Section~\ref{sub:fsvt}.  Moreover, this $\mathit{K}$-polynomial and
multidegree are known \cite{grobGeom}: they are, respectively, the
\emph{double Grothen\-dieck polynomial} and \emph{double Schubert
polynomial} associated by Lascoux and Sch\"utzenberger \cite{LS82}
to~$\pi$; see Section~\ref{sub:flaggings} for our conventions on
Schubert and Grothendieck polynomials in this paper. 
The previous theorem implies the following.

\begin{Corollary*}
Let $\pi \in S_n$ be a vexillary permutation.  
The double Schubert polynomial can be expressed as the positive sum
\begin{eqnarray*}
  \SS_\pi(\xx,\yy) &=& \ \sum_{\tau\in \FT(\pi)\,} \prod_{e\in \tau} \ 
  (x_{\val(e)} - y_{\val(e) + j(e)})
\end{eqnarray*}
over the set $\FT(\pi)$ of flagged tableaux associated to~$\pi$, 
all of which have shape~$\lambda(\pi)$.  
Each product is over the entries $e$ of $\tau$, whose numerical
values are denoted $\val(e)$, and
where $j(e) = c(e) - r(e)$ is the difference of the row and column
indices. 

The double Grothendieck
polynomial $\GG_\pi(\xx,\yy)$ can be expressed as the % noncanceling
sum
\begin{eqnarray*}
  \GG_\pi(\xx,\yy) 
  &=& \sum_{\tau\in \FST(\pi)\,} (-1)^{|\tau|-|\lambda|}
  \prod_{e\in \tau} \Big(1 - \frac{x_{\val(e)}}{y_{\val(e)+j(e)}} \Big)  
\end{eqnarray*}
over the set $\FST(\pi)$ of flagged set-valued tableaux associated
to~$\pi$.
The sign $(-1)^{|\tau|-|\lambda|}$ alternates with the number of
``excess'' entries in the set-valued tableau.
\end{Corollary*}

This result appears as part of Theorem~\ref{t:FSVT}.
% (before which we lay out our conventions for Grothendieck
% polynomials).
It was already known in the case of single Schubert polynomials
\cite{Wachs} and in the case of Grothendieck polynomials for
Grassmannian permutations \cite{Buch}.
% in the antidiagonal case ... recover the pipe dream formulae of
% \cite{FKyangBax,BB} ...

Since the paper \cite{grobGeom} already provided a geometric
explanation for a combinatorial formula for the Grothendieck and
Schubert polynomials of any vexillary permutation---indeed, of any
permutation---the reader may wonder why we have provided another
one. There are three reasons. The primary reason is to show that the
combinatorial trick we used in \cite{subword}, vertex-decomposability,
can have a transparent geometric origin. Another is to directly
connect the Gr\"obner geometry to Young tableaux, rather than to the
less familiar pipe dreams. Finally, the two formulae themselves are
very different, as demonstrated in Example~\ref{exa:hard_bij}.

One of the satisfying aspects of the degenerations in this paper is
that they stay within the class of (vexillary) matrix Schubert varieties.
While one \emph{can} view the antidiagonal degeneration used in
\cite{grobGeom} also in terms of geometric vertex decompositions, 
we don't know what the cone and projection pieces look like along the way;
they are no longer matrix Schubert varieties.

%%%%%%%%%%%%%%%%%%%%%%%%%%%%%%%%%%%%%%%%%%%%%%%%%%%%%%%%%%%%%%%%%%%%%%
\section{Geometric vertex decompositions}%%%%%%%%%%%%%%%%%%%%%%%%%%%%%
%%%%%%%%%%%%%%%%%%%%%%%%%%%%%%%%%%%%%%%%%%%%%%%%%%%%%%%%%%%%%%%%%%%%%%

The main results on geometric vertex decompositions are most easily
stated in algebraic language, so we do this first, in
Theorem~\ref{thm:algGVD}.  Then, for ease of future reference, we make
explicit in Section~\ref{sub:geom} the geometric interpretation of
Theorem~\ref{thm:algGVD}.  Our choice of geometric language makes it
clear that the description of geometric vertex decomposition already
given in Section~\ref{ssec:gvd} does not really depend on the
hyperplane~$H$.  We close Section~\ref{sub:geom} with a useful
technique for working with (repeated) geometric vertex decompositions
of reduced schemes.

\subsection{Algebraic aspects}

Let $R = \kk[x_1,\ldots,x_n,y]$ be a polynomial ring in $n+1$
variables over an arbitrary field~$\kk$.  We shall be dealing with
Gr\"obner bases, for which our basic reference is
\cite[Chapter~15]{eisenbud}.  Define the \dfn{initial $y$-form}
$\init_y p$ of a polynomial $p\in R$ to be the sum of all terms of~$p$
having the highest power of~$y$.  Fix a term order~$\Le$ on~$R$ such
that the initial term $\init\,p$ of any polynomial~$p$ is a term in
the initial $y$-form: $\init\,p = \init(\init_y p)$.  Let $\init\,J$
denote the initial ideal of $J$, generated by the initial terms
$\init\,f$ of all $f\in J$, and let $\init_y J$ denote the ideal
generated by the initial $y$-forms of the elements of~$J$.  Thus
$\init\,\<g_1,\ldots,g_r\> \supseteq \<\init\,g_1, \ldots,
\init\,g_r\>$, with equality if and only if $\{g_1,\ldots,g_r\}$ is a
Gr\"obner basis under~$\Le$.  We automatically have $\init(\init_y J)
= \init\,J$, by our condition on~$\Le$.

We say that $I$ is \dfn{homogeneous} if it is $\ZZ$-graded for the grading
on~$R$ in which all of the variables have degree~$1$.  When $I$ is
homogeneous, it has a Hilbert series $h_{R/I}(s) = \sum_{k\in\ZZ}
\dim(R/I)_k \cdot s^k$.  When $h$ and $h'$ are two power series in the
variable~$s$ with integer coefficients, we write $h \leq h'$ if for
all~$k$, the coefficient on~$s^k$ in~$h$ is less than or equal to the
coefficient on~$s^k$ in~$h'$.

\begin{Theorem}\label{thm:algGVD}
  Let $I$ be an ideal in $R$, and $\{y^{d_i} q_i + r_i\}_{i=1}^m$
  a Gr\"obner basis for it, where $y^{d_i} q_i$ is the initial
  $y$-form of $y^{d_i} q_i + r_i$ and $y$ does not divide~$q_i$.  Then
  the following statements hold for
  \[
    \qquad I' = \<y^{d_i} q_i \mid i = 1,\ldots,m\>, \quad C = \<q_i
    \mid i = 1,\ldots,m\>, \quad P = \<q_i \mid d_i = 0\> + \<y\>.
  \]
  \begin{alphlist}
  \item
  The given generating sets of~$I'$, $C$, and~$P$ are Gr\"obner bases,
  and $I' = \init_y I$.
  \item
  If $\max_i d_i = 1$, then $I' = C\cap P$ and $\init(C+P) = \init\,C
  + \init\,P$.
  \item
% $I' \subseteq C \cap P$, and
  $\sqrt{I'} = \sqrt{C}\cap\sqrt{P}$.
  \item
  $C = (I' : y^\infty)$, which is by definition the ideal $\{f \in R
  \mid y^j f \in I'$ for some~$j \geq 0\}$.
  \item
  If $I$ is homogeneous, then $h_{R/I}(s) \geq h_{R/P}(s) + s\cdot
  h_{R/C}(s)$, with equality if and only if $I' = C\cap P$.
  \end{alphlist}
When $I' = C\cap P$, we will call this decomposition of $I'$ a
\dfn{geometric vertex decomposition of~$I$}.
\end{Theorem}

\begin{proof}
(a)
These all follow from
% \cite[Proposition~15.29 and Proposition~15.31]{eisenbud}.
\cite[Section~15.10.4 and~15.10.5]{eisenbud}.
  
(b)
% {\bf How do we use $yC$? Why confuse things?}
% EM: I hope the rephrasing makes this clear.
The hypothesis means that $I' = yC + \<q_i \mid d_i = 0\>$.  As the
generators of~$C$ do not involve the variable~$y$, we have $yC = C
\cap \<y\>$.  Hence, using the modular law for ideals (that is,
$\mathfrak{c \cap (a+b) = c \cap a + c \cap b}$ if $\mathfrak{c
\supseteq b}$), we can conclude that $I' = C \cap \<y\> + C \cap \<q_i
\mid d_i = 0\> = C \cap P$.  For the equality involving initial
ideals, observe that $C + P = C + \<y\>$, so that $\init(C+P) =
\init\,C + \<y\>$, and use part~(a).

(c) Let $\tilde I' = \<y q_i \mid d_i \geq 1\> + \<q_i \mid d_i =
0\>$.  Then $\sqrt{I'} \supseteq \tilde I' \supseteq I'$, so by the
Nullstellen\-satz all three ideals have the same vanishing set.
Part~(b) implies $\tilde I' = C \cap P$.

(d) This is elementary, using the fact that $I'$ is homogeneous for
the $\ZZ$-grading in which $y$ has degree~$1$ and all other variables
have degree zero.

(e) It suffices to show that the sum of series on the right hand side
of the inequality is the Hilbert series of $R/(C \cap P)$, because
(i)~$C \cap P \supseteq I'$ by our proof of~(c), and (ii)~the
quotients $R/I$ and $R/I'$ have the same Hilbert series.  To complete
the proof, use the exact sequence
\[
  0 \to R/(C+P) \to R/C \oplus R/P \to R/(C\cap P) \to 0
\]
of $R$-modules.  It implies that
\[
  h_{R/(C\cap P)} = h_{R/P} + (h_{R/C} - h_{R/(C+P)}).
\]
The equality $C+P = C + \<y\>$ yields $R/(C+P) = (R/C)/y(R/C)$.
Therefore $h_{R/(C+P)}(s) = (1-s)h_{R/C}(s)$, because the generators
of~$C$ do not involve~$y$.%
\end{proof}

This inequality (e) is used in \cite{KMN} to study schemes whose Hilbert
functions are smallest, in various senses.

\subsection{Geometric aspects}\label{sub:geom}

\newcommand\bP{{\mathbb P}}
\newcommand\tP{\mathit{Bl}_L\ol V}

While the next theorem essentially recapitulates
Theorem~\ref{thm:algGVD} in geometric language, it is not a verbatim
translation.  For example, we do not assume that coordinates
$x_1,\ldots,x_n,y$ have been given.
% we do not attempt to translate part~(b) of
% Theorem~\ref{thm:algGVD}.
One of the purposes of Theorem~\ref{thm:geoGVD} is to describe the
flat limit~$X'$ using schemes naturally determined by the subscheme $X
\subseteq V$ and the choice of the line~$L$, namely $\Pi$
and~$\Lambda$, at least in the case where we have a geometric vertex
decomposition.

Let $V$ be a vector space over a field~$\kk$, viewed as a scheme
over~$\kk$, and suppose that a $1$-dimensional subspace~$L$ of~$V$ has
been given.  The projectivization $\ol V = \bP(V \oplus \kk)$, which
we view as the projective completion $V \cup \bP V$ of~$V$, has a
point $\bP L \in \bP V$.  Denote by~$\tP$ the blow up of $\bP(V \oplus
\kk)$ at the point~$\bP L$.  The exceptional divisor is naturally
identified with the projective completion of the quotient vector
space~$V/L$,
% {\bf is this true? or is the natural thing $\bP(V/L \oplus L)$,
%  and I should only be doing this in the homogeneous case? Nah.}
and in particular contains a copy of~$V/L$.

For each choice of a codimension~$1$ subspace $H \subseteq V$
complementary to~$L$, there is an action $t \cdot (\vec h,l) = (\vec
h,tl)$ of\/~$\kk^\times$ on~$V$, which we call \dfn{scaling~$L$ and
fixing~$H$}.
% Observe that $H$ maps isomorphically to~$V/L$ under the projection
% of~$V$ modulo~$L$.

\begin{Theorem}\label{thm:geoGVD}
  Let $X$ be a closed subscheme of\/~$V\!$ and $L$ a $1$-dimensional
  subspace of\/~$V$.  Denote by $\Pi$ the scheme-theoretic closure of
  the image of\/~$X$ in~$V/L$, and by~$\ol X$ the closure of\/~$X$
  in~$\tP$.  Set $\Lambda = \ol X \cap V/L$, where the intersection of
  schemes takes place in~$\tP$.
  \begin{alphlist}
  \item
If~$H$ is a hyperplane complementary to $L$ in~$V$, and we
identify~$H$ with~$V/L$, then the flat limit $X' := \lim_{t\to 0}
t\cdot X$ under scaling~$L$ and fixing~$H$ satisfies
\[
  X' \supseteq (\Pi \times \{0\}) \cup_{\Lambda \times \{0\}} (\Lambda
  \times L).
\]
\item
The scheme-theoretic containment in part~(a) is an equality as sets.
\item
If (the ideal of) $X$ is homogeneous, then the same holds for~$\Pi$ as
well as~$\Lambda$, and we derive an inequality on Hilbert series of
subschemes of\/~$V$:
\[
  h_X(s) \geq (1-s) h_{\Pi\times L}(s) + s h_{\Lambda\times L}(s).
\]
\item
Parts~(a) and~(c) both become equalities if\/ the flat limit $X'$ is
reduced.
\end{alphlist}
\end{Theorem}
\begin{proof}
% The description of $\tP$ is standard.
% Note that the action defined of $\kk^\times$ extends to $\tP$, and
% the fixed points are $\bP(H\oplus \kk)$ and the exceptional
% divisor. {\bf I don't want to get deep into TV stuff to identify the
% open set $V \cup V/L$ with $H\times (L\cup\infty)$... what's a
% simple way to say this?}
% EM: I don't see a need to say anything about it at all; the reader
% wishing to make the connection to the $H \times L\cup\infty$
% language we used in the Introduction ought to be able to do it
% themselves.  The only thing we need to derive the statements from
% Theorem~\ref{thm:algGVD}.  
Pick coordinates $x_1,\ldots,x_n$ on $H$ and a coordinate $y$ on $L$,
and choose a term order on $R = \kk[x_1,\ldots,x_n,y]$ such that
$\init(\init_y\, p) = \init\, p$ for all polynomials~$p$.
% EM: the next (now %-commented out) phrase is unnecessary, and in any
% case is unclear, because the order of the variables has to be specified: 
% ; for example a lexicographic term order.
Let $I \subseteq R$ be the ideal defining $X$.  Then there exists a
Gr\"obner basis for~$I$ with respect to this term order, and we can
apply Theorem~\ref{thm:algGVD} to study the associated ideals $I'$,
$C$, and~$P$.  The ideal $P$ cuts out the projection $\Pi \subseteq
H$, while $C$ cuts out the subscheme $\Lambda \times L \subseteq V$.
% Parts~(c) and~(e) of Theorem~\ref{thm:algGVD} now imply all of
% our~claims.
Our claims therefore follow from Theorem~\ref{thm:algGVD}.%
\end{proof}

Though it is very important in this paper, we did not see how to
state Theorem~\ref{thm:algGVD} part~(b) in a particularly geometric way.  
We expressed the Hilbert series in Theorem~\ref{thm:geoGVD}(c)
in terms of $X$, $\Pi\times L$, and $\Lambda\times L$ because all
three occur in the same vector space (once $H$ has been chosen).

% \subsection{Stability of geometric vertex decompositions}

The property of being a geometric vertex decomposition is preserved
under (further) degeneration, so long as the schemes stay reduced:
%The example we have in mind for the decomposition of~$Y$ in what
%follows is a geometric vertex decomposition of another scheme in
%$H\times L\times M$.

\begin{Proposition}\label{prop:inductivegvd}
  Let $X\supseteq D$ be two reduced closed subschemes of $H\times L$.
  Let $M$ be a $1$-dimensional vector space, and define $Y\subseteq
  H\times L\times M$ as
  \[
    Y = (X \times \{0\}) \cup_{D\times \{0\}} (D\times M).
  \]
  Use the action of~$\kk^\times$ on~$H \times L$ scaling~$L$ and
  fixing~$H$ to define two new flat limits $X'$ and~$D'$.  Assume that
  these are again reduced,
% EM: is reduced equivalent to being a geometric vertex decomposition?
  so they are geometric vertex decompositions.
% Then if we use the decomposition into $H \times M$ and $L$, and the
% rescaling action on $L$ to define a new flat limit $Y'$ by rescaling
% $Y$, then $Y'$ is again reduced, i.e.\ a geometric vertex
% decomposition.
  Then the flat limit $Y'$ of~$Y$ under the action of~$\kk^\times$ on
  $H \times M \times L$ scaling~$L$ and fixing $H\times M$ is again
  reduced, and hence it is a geometric vertex decomposition.%
\end{Proposition}

\begin{proof}
  By Theorem~\ref{thm:geoGVD}, as subschemes of $H\times L$ we have
  \begin{eqnarray*}
     X' = (\Pi\times \{0\}) \cup (\Lambda \times L) &\hbox{and}& D' =
     (\Sigma\times \{0\}) \cup (\Gamma \times L).
  \end{eqnarray*}
  Hence the projection and cone parts of $Y$' are, as subsets of~$H
  \times M$,
  \begin{eqnarray*}
    (\Pi\times \{0\}) \cup_{\Sigma\times \{0\}} (\Sigma \times M)
    &\hbox{and}& (\Lambda \times \{0\}) \cup_{\Gamma \times \{0\}}
    (\Gamma \times M),
%   \qquad \subseteq H\times M,
  \end{eqnarray*}
  respectively.  Again by Theorem~\ref{thm:geoGVD} we get
  \[
    Y' \supseteq
    (\Pi\times \{0\}\times \{0\}) 
    \cup (\Sigma \times \{0\} \times M)
    \cup (\Lambda \times L \times \{0\})
    \cup (\Gamma \times L \times M).
  \]
  Our goal is to prove the above containment to be an equality of
  schemes.

  Rearranging, we see that
  \begin{eqnarray*}
    Y' &\supseteq&
    (\Pi\times \{0\}\times \{0\}) 
    \cup (\Lambda \times L \times \{0\})
    \cup (\Sigma \times \{0\} \times M)
    \cup (\Gamma \times L \times M) \\
    &=& (X' \times \{0\}) \cup_{D'\times \{0\}} (D' \times M).
  \end{eqnarray*}
  Since $X'$ and $D'$ had geometric vertex decompositions, the right
  hand side has the same Hilbert series as $(X \times \{0\})
  \cup_{D\times \{0\}} (D\times M) = Y$, which matches that of $Y'$,
  the left hand side. Therefore the containment is an equality.%
% \comment{this used to say ``inequality'', which EM presumed to be a
% typo}
\end{proof}

\subsection{Cohomological aspects}
\newcommand\Sym{{\rm Sym}}
\newcommand\codim{{\rm codim}}

Even if the degeneration by rescaling an axis is \emph{not}\/ a
geometric vertex decomposition, the limit $X'$ can still be analyzed
enough to give a useful positivity statement about (multi)degrees.
Moreover, it is always an equality, not just an inequality like
Theorem~\ref{thm:geoGVD} part~(c).

For $X$ a $T$-invariant subscheme of a vector space~$V$ carrying an
action of a torus~$T$, there is an associated \dfn{multidegree}
$\deg_{V\!} X$ living in the symmetric algebra $\Sym(T^*)$, where
$T^*$ is the weight lattice of~$T$.  Our general reference for
multidegrees is \mbox{\cite[Chapter~8]{cca}}, though this particular
algorithm appears in \cite{Jo97}.

\begin{Proposition} \label{p:multideg}
Three axioms suffice to characterize the assignment of multidegrees:
\begin{alphlist}
\item $\deg_{\{\vec 0\}} \{\vec 0\} = 1$ for the zero vector space
  $\vec 0$ and the trivial torus action.
\item If the components of $X$ of top dimension are $X_1,\ldots,X_k$,
  occurring in $X$ with multiplicities $m_1,\ldots,m_k$, then $\deg_V
  X = \sum_i m_i \deg_{V\!} X_i$.
\item If $H$ is a $T$-invariant hyperplane in $V$, and $X$ is an
  irreducible variety, then
  \begin{romanlist}
  \item if $X \not\subseteq H$, then $\deg_{V\!} X = \deg_H(X\cap H)$.
  \item if $X \subseteq H$, then $\deg_{V\!} X = \wt(V/H)\deg_H X$,
    where $\wt(V/H) \in \Sym^1(T^*)$ is the weight of the $T$-action
    on the $1$-dimensional representation~$V/H$.
  \end{romanlist}
\end{alphlist}
\end{Proposition}
\begin{proof}
That these properties are satisfied by multidegrees follows from
\cite[Section~8.5 and Exercise~8.12]{cca}.  These properties determine
a unique assignment by induction on the dimension of~$V$: the base
case is the first axiom; the second axiom reduces the calculation of
multidegrees from schemes to varieties;
% (reduced irreducible schemes);
and the third axiom brings the dimension down by~1 for varieties.%
\end{proof}

If $T$ is $1$-dimensional, then $\Sym^{\codim X}(T^*) \iso \integers$,
so $\deg_{V\!} X \in \Sym^{\codim X}(T^*)$ is specified by a number.
If in addition the torus action is by global rescaling on $V$, so $X$
is an affine cone, then this number is the usual degree of the
corresponding projective variety. This cohomological interpretation is
related to the general fact that $\Sym(T^*)$ is the equivariant
cohomology ring, and equivariant Chow ring, of $V$.  The theorem to
follow is already interesting as a statement about usual degrees.

In an unfortunate collision of terminology, the \dfn{degree} of a
morphism%
        \footnote{Actually, the two are related. If $W\leq V$ is a
        generic subspace of codimension $\dim X$, then the degree of
        the morphism $X \to V/W$ is the usual degree of the affine
        cone $X$. Such genericity is unavailable for most bigger torus
        actions, and one may view Theorem \ref{thm:deg} as a version
        of this statement for multidegrees.}
$X \to Y$ of reduced irreducible schemes over a field~$\kk$ is defined
to be the degree of the extension $\kk(X) \supseteq \kk(Y)$ of their
fraction fields (if this extension is finite) and zero otherwise.
When $\kk$ is algebraically closed, the degree of $X \to Y$ is simply
the cardinality of a generic fiber, if this number is finite.%
\vspace{-3.4pt}

\goodbreak

\begin{Theorem} \label{thm:deg}
  Let $X,X' \subseteq H\times L$ and $\Pi,\Lambda \subseteq H$ be as in
  Theorem~\ref{thm:geoGVD}.  Assume that $X$ is reduced, irreducible,
  and invariant under the action of a torus $T$ on $V$, so its
  projection $\Pi \subset H$ is, too.
% Assume the following equivalent conditions: $X$ is not invariant
% under translation by $L$, $\Pi \supsetneq \Lambda$, the map
% $X\to\Pi$ has generically finite fibers.  The map $X\to \Pi$
% therefore has a well-defined degree~$d$.
  Let $d$ be the degree of the projection morphism $X \to \Pi$.  Then
  \[
  \deg_{V\!} X = d \cdot \deg_V \Pi + \deg_V (\Lambda \times L).
  \]
  Moreover, $\deg_V \Pi = \deg_V (\Pi \times L) \cdot \wt(V/H)$,
  where $\wt(V/H) \in \Sym^1(T^*)$ is the weight of the action of $T$
  on the one-dimensional representation $V/H$.
\end{Theorem}

\begin{proof}
% Algebraically, the degree of the projection map $X \to \Pi$ is the
% rank of the algebra $Fun(X)$ as a module over the domain $Fun(\Pi)$.
% This rank doesn't change when we replace $X$ by $X'$, which also
% maps to $\Pi$.
% 
% \renewcommand\tP{\widetilde{P}}
% Recall the ideals $I' = \<y^{d_i} q_i \mid i = 1,\ldots,m\>, C =
% \<q_i \mid i = 1,\ldots,m\>, P = \<q_i \mid d_i = 0\> + \<y\>$ from
% Theorem \ref{thm:algGVD}.  Let $k = \max_i d_i$. Let $\tP = I' +
% \<y^k\>$.
% 
% Then we claim that $I' = \tP\ \cap\ C$.  Both terms on the RHS
% contain $I'$, so we need to show $\supseteq$.  Let $p\in \tP$, so $p
% = y^k a + b$ where $b\in I'$.  Therefore $y^k a \in I' \leq C$,
% hence $a\in C$, since no generator of $C$ involves $y$. If $a =
% \sum_i m_i q_i$, then $y^k a = \sum_i (m_i y^{k-d_i}) (y^{d_i} q_i)
% \in I'$. So $p\in I'$.
% 
% This gives a relation on Hilbert series:
% $$
%   h_{R/I} = h_{R/I'} = h_{R/\tP} + h_{R/C} - h_{R/(\tP + C)}.
% $$
% The scheme defined by $\tP + C$ is supported on the intersection
% $\Lambda \times \{0\}$, which by assumption has lower dimension than
% $\Pi \times \{0\}$.  Hence we can neglect it when computing degrees:
% $$
%   \deg R/I = \deg R/\tP + \deg R/C.
% $$
% It remains to show that $\deg R/\tP = d \deg \Pi$.
% 
% Plainly $\sqrt{\tP} \geq P \geq \tP$, so the support of $P$ and
% $\tP$ is the same, $\Pi$. {\bf now what to quote?}
  Let us assume that the morphism $X \to \Pi$ has degree $d > 0$, for
  otherwise $X = \Pi \times L = X'$, and the result is trivial.  The
  stability of $\deg$ under Gr\"obner degeneration
  \cite[Corollary~8.47]{cca} implies that $\deg(X) = \deg(X')$.
  Additivity of $\deg$ on unions \cite[Theorem~8.53]{cca} therefore
  implies, by Theorem~\ref{thm:geoGVD}, that $\deg(X) = \deg(\Lambda
  \times L) + \deg(\Pi')$, where $\Pi'$ is the (possibly nonreduced)
  component of~$X'$ supported on~$\Pi$.  It remains only to show that
  $\deg(\Pi') = d \cdot \deg(\Pi)$.  This will follow from the
  additivity of~$\deg$ once we show that~$\Pi'$---or
  equivalently,~$X'$---has multiplicity~$d$ along~$\Pi$.
  
  Let $K = \kk(\Pi)$, the fraction field of~$\Pi$.  Given a module~$M$
  over the coordinate ring~$\kk[\Pi]$ of~$\Pi$, denote by~$M_\Pi$ the
  localization $M \otimes_{\kk[\Pi]} K$ at the generic point of~$\Pi$.
  The multiplicity of~$X'$ along~$\Pi$ is the dimension
  of~$\kk[X']_\Pi$ as a vector space over~$K$.  Now consider the flat
  degeneration $X \rightsquigarrow X'$ as a family over the line with
  coordinate ring~$\kk[t]$, and let $M$ be the coordinate ring of the
  total space of this family.  Thus $M$ is flat as a module
  over~$\kk[t]$; equivalently, $M$ is torsion-free over~$\kk[t]$.  To
  prove the result it is enough to show that the localization $M_\Pi$ is
  flat over~$K[t]$, since taking $K$-vector space dimensions of the
  fibers over $t = 1$ and $t = 0$ yields the degree~$d$ of $X \to \Pi$
  and the multiplicity of~$X'$ along~$\Pi$, respectively.
  
  We know that $M$ is flat over~$\kk[t]$, that $\kk[\Pi][t]$ is flat
  over~$\kk[t]$, and that $K[t]$ is flat over~$\kk[\Pi][t]$.  A routine
  calculation therefore shows that $M_\Pi$ is flat over~$\kk[t]$.
  On the other hand, $M$ is the coordinate ring of the
  total space of a (partial) Gr\"obner degeneration, so $M
  \otimes_{\kk[t]} \kk[t,t^{-1}]$ is the coordinate ring
  of a family over $\kk^\times$ that is isomorphic to the
  trivial family $\kk^\times \times X$.  Moreover, since
  the rescaling in our Gr\"obner degeneration commutes
  with the projection to~$\Pi$, the coordinate ring $M
  \otimes_{\kk[t]} \kk[t,t^{-1}]$ is free over
  $\kk[\Pi][t,t^{-1}]$.  It follows that $M_\Pi \otimes_{K[t]} K[t,t^{-1}]$
  is free over~$K[t,t^{-1}]$.  Consequently, the $K[t]$-torsion
  submodule of~$M_\Pi$ is supported at $t = 0$.  The $K[t]$-torsion
  of~$M_\Pi$ is thus also $\kk[t]$-torsion, each element therein being
  annihilated by some power of~$t$.  Since $M_\Pi$ has no
  $\kk[t]$-torsion, we conclude $M_\Pi$ is torsion-free, and hence flat,
  over~$K[t]$.

% The second equation follows from properties (3a) and (3b) of
% multidegrees. 
  The claim about $\deg_V\Pi$ in the final sentence follows from
  Proposition~\ref{p:multideg}.%
\end{proof}

We remark that in the setup of Theorem~\ref{thm:deg}, 
$d=1$ is necessary in order to have a geometric 
vertex decomposition. It is not sufficient, as seen in 
Example \ref{ex:hyperbolaplust}.

This theorem is applied in \cite{Ksubword} to affine
patches on Schubert varieties.

For a subtorus $S \subseteq T$, there are two different multidegrees
\mbox{$\deg_{V\!}^{T\!} X \in \Sym(T^*)$} and $\deg_V^S X \in
\Sym(S^*)$ one could assign to a $T$-invariant (hence also
$S$-invariant) subscheme $X\subseteq V$.  They are related by the
natural map $\Sym(T^*)\to\Sym(S^*)$ induced from the restriction of
characters $T^*\to S^*$.

The map $\Sym(T^*)\to\Sym(S^*)$ can be thought of as specialization of
polynomials upon the imposition of linear conditions on the variables.
The following consequence of Theorem~\ref{thm:deg} gives a criterion
for one polynomial to be the specialization of another.  In it, we do
not assume any Gr\"obner properties of $\{y^{d_i}q_i + r_i\}$; hence
the ideal $J$ here is contained in the ideal~$C$ from
Theorem~\ref{thm:algGVD}, perhaps properly.

\begin{Corollary}
  Let $X \subseteq H\times L$, where $H$ has coordinates
  $x_1,\ldots,x_n$ and $L$ has coordinate~$y$.  Assume that $H$
  and~$L$ are representations of a torus~$T$, and $X$ is a
  $T$-invariant subvariety.  Let $w \in T^* = \Sym^1(T^*)$ be the
  weight of~$T$ on~$L$, and $S \subseteq T$ the stabilizer of~$L$, so
  the map $\Sym(T^*)\to \Sym(S^*)$ takes $p\mapsto p|_{w=0}$.

  Let the ideal $I$ defining $X$ be generated by $\{y^{d_i} q_i + r_i
  \mid i = 1,\ldots,m\}$, where $y^{d_i} q_i$ is the initial $y$-form
  of $y^{d_i} q_i + r_i$ and $y$ does not divide~$q_i$.  Let $J =
  \<q_1,\ldots,q_m\>$,
% (Since we have not assumed any Gr\"obner properties of $\{y^{d_i}q_i
% + r_i\}$, this $J$ is contained in the $C$ from
% Theorem~\ref{thm:algGVD}, perhaps properly.)
  and define $\Theta \subseteq H$ to be the zero scheme of~$J$.  If we
  know that
  \begin{itemize}
  \item $\Theta$ has only one component of dimension $\dim X - 1$,
  \item that component is generically reduced, and
  \item $X$ is not contained in a union of finitely many translates of $H$,
  \end{itemize}
  then
  $$ (\deg_{V\!} X)|_{w=0} = (\deg_H \Theta)|_{w=0}. $$
\end{Corollary}

\begin{proof}
  Let $\Theta'$ be the reduced variety underlying the $(\dim
  X-1)$-dimensional component of $\Theta$. By the conditions on
  $\Theta$, $\deg_H \Theta = \deg_H \Theta'$.

  Let $\Lambda,\Pi,d$ be as in Theorem \ref{thm:deg}, 
  whose conclusion specializes to
  $$ 
  (\deg_{V\!} X)|_{w=0} = (\deg_H \Lambda)|_{w=0}
  $$
  since setting $w=0$ kills the contribution from $\Pi$.

  Let $C$ be the ideal defining $\Lambda$. Then $C$ contains~$J$, so
  $\Lambda \subseteq \Theta$.  The condition that $X$ not be contained
  in a union of finitely many translates of $H$ says that $\Lambda$ is
  nonempty, and hence has dimension $\dim X - 1$.

  Hence $\Lambda \supseteq \Theta'$, and being trapped between two schemes
  with the same multidegree, $\deg_H \Lambda = \deg_H \Theta$.
\end{proof}

Essentially, this corollary replaces the difficulty of showing that 
a basis for an ideal is Gr\"obner with the difficulty of showing that
$\Theta$ has only one big component. It is used in precisely this form
in \cite{KZJ}.

%%%%%%%%%%%%%%%%%%%%%%%%%%%%%%%%%%%%%%%%%%%%%%%%%%%%%%%%%%%%%%%%%%%%%%
\section{Vexillary matrix Schubert varieties}\label{sec:gvdvmSv}%%%%%%
%%%%%%%%%%%%%%%%%%%%%%%%%%%%%%%%%%%%%%%%%%%%%%%%%%%%%%%%%%%%%%%%%%%%%%

\subsection{Matrix Schubert varieties} \label{sub:matrixSchub}

In this subsection we review some definitions and results of Fulton on
determinantal ideals \cite{fulton:92}; an exposition of this material
can be found in \cite[Chapter~15]{cca}.

Let $M_n$ be the variety of $n\times n$ matrices over $\kk$, with
coordinate ring $\kk[\zz]$ in indeterminates $\{z_{ij}\}_{i,j=1}^{n}$.
We will let $\zz$ denote the generic matrix of variables $(z_{ij})$
and let $\zz_{p\times q}$ denote the northwest $p\times q$ submatrix
of~$\zz$.  More generally, if $Z$ is any rectangular array of objects,
let $Z_{p \times q}$ denote the northwest $p \times q$ subarray.  In
particular, identifying $\pi \in S_n$ with the square array having
blank boxes in all locations except at $(i,\pi(i))$ for $i =
1,\ldots,n$, where we place \dfn{dots}, we define the rank $\rr pq\pi$
to be the number of dots in the subarray~$\pi_{p \times q}$.  This
yields the $n \times n$ \dfn{rank array} $\rrr\pi = (\rr
pq\pi)_{p,q=1}^n$.  We can recover the dot-matrix for~$\pi$
from~$\rrr\pi$ by placing a dot at $(p,q)$ whenever $\rr pq\pi -
\rr{p-1,}{q}\pi - \rr{p,}{q-1}\pi + \rr{p-1,}{q-1}\pi = 1$ rather than $0$ (by
convention, we set $\rr 0*\pi = \rr *0\pi = 0$).

% \comment{need to fix the $r_{??}(?)$ notation here:}
% Define a \emph{rank matrix} as a matrix $r$ of natural numbers,
% satisfying three conditions:
% \begin{itemize}
% \item $r_{pq} \leq p,q$
% \item $r_{p+1,q}, r_{p,q+1} \in [r_{pq},1+r_{pq}]$
% \item $r_{p+1,q+1} + r_{pq} \geq r_{p+1,q} + r_{p,q+1}.$
% \end{itemize}
% \comment{I wonder if we need this way of talking about rank
% functions. Of course we could just directly define $r_\pi$ as the
% rank function for $\pi$ obtained by computing the number of dots in
% northwest rectangles etc.} To each such $r$ we can associate a $0,1$
% matrix $\pi$ by $\pi_{pq} = r_{pq} - r_{p-1,q} - r_{p,q-1} +
% r_{p-1,q-1}$ (taking $r_{0*} = r_{*0} = 0$). Then $\pi$ is a partial
% permutation matrix, and the map $r \mapsto \pi$ provides a bijective
% correspondence between rank matrices and partial permutation
% matrices. 
% If we let $G(\pi)$ denote the permutation matrix associated to
% $\pi$, and call the nonzero entries of $G(\pi)$ its \emph{dots}, then
% the inverse map $\pi\mapsto r^{(\pi)}$ is given by setting
% $r^{(\pi)}_{ij}$ equal to the number of dots appearing the in the
% northwest $i\times j$ rectangle of $G(\pi)$.

For $\pi \in S_n$, the \dfn{Schubert determinantal ideal} $I_\pi
\subseteq \kk[\zz]$ is generated by all minors in $\zz_{p\times q}$ of
size $1+\rr pq\pi$ for all $p$ and~$q$.  It was proven in \cite{fulton:92}
to be prime. The \dfn{matrix Schubert
variety} $\ol X_\pi$ is the subvariety of~$M_n$ cut out by~$I_\pi$; 
thus $\ol X_\pi$ consists of all matrices $Z\in M_{n}$
such that $\mathrm{rank}(Z_{p\times q})\leq \rr pq\pi$ for all $p$
and~$q$.

In fact, the ideal $I_\pi$ is generated by a smaller subset of these
determinants.  This subset is described in terms of the \dfn{diagram}
\[
  D(\pi)=\big\{(p,q) \in \{1,\ldots,n\} \times \{1,\ldots,n\} \mid
  \pi(p)>q \hbox{ and } \pi^{-1}(q)>p\big\}.
\]
of~$\pi$.  Pictorially, if we draw a ``hook'' consisting of lines
going east and south from each dot, then $D(\pi)$ consists of the
squares not in the hook of any dot.
% The diagram appears as a collection of connected components
The \dfn{essential set} is the set of southeast corners of the
connected components of the diagram:
\[
  \Ess(\pi)=\{(p,q)\in D(\pi) \mid (p+1,q),\ (p,q+1)\not\in D(\pi)\}.
\]
Then with these definitions, we call a generator of~$I_\pi$
\dfn{essential} if it arises as a minor of size $1+\rr pq\pi$ in
$\zz_{p\times q}$ where $(p,q)\in \Ess(\pi)$.  The prime ideal~$I_\pi$
is generated by its set of essential minors.

% We are interested in diagonal term orders on $\kk[\zz]$, which by
% definition is pick off from each minor its diagonal term.  Many
% examples of these term orders exist, including, e.g., reverse
% lexicographic order on the variables $\zz$ as read from right to
% left and bottom-up.

\begin{Example}\label{exa:3142}
The dot-matrix for~$\pi= \left(\begin{array}{ccccc} 1 & 2 & 3 & 4 &
5\\ 4 & 1 & 3 & 2 & 5 \end{array}\right)\in S_{5}$ and the diagram
$D(\pi)$ are combined below:
% \[
% \begin{picture}(200,100)
% \put(-30,40){$G(\pi),\ D(\pi)=$}
% \put(50,0){\makebox[0pt][l]{\framebox(100,100)}}
% \thicklines
% \put(60,70){\circle*{4}}
% \put(60,70){\line(1,0){90}}
% \put(60,70){\line(0,-1){70}}
% 
% \put(80,30){\circle*{4}}
% \put(80,30){\line(1,0){70}}
% \put(80,30){\line(0,-1){30}}
% \put(120,90){\circle*{4}}
% \put(120,90){\line(1,0){30}}
% \put(120,90){\line(0,-1){70}}
% \put(100,50){\circle*{4}}
% \put(100,50){\line(1,0){50}}
% \put(100,50){\line(0,-1){50}}
% \put(140,10){\circle*{4}}
% \put(140,10){\line(1,0){10}}
% \put(140,10){\line(0,-1){10}}
% \thinlines
% \put(70,40){\makebox[0pt][l]{\framebox(20,20)}}
% \thicklines
% \put(120,30){\line(1,0){10}}
% \put(120,30){\line(0,-1){30}}
% \thinlines
% \put(50,80){\line(1,0){60}}
% \put(90,80){\line(0,1){20}}
% \put(70,80){\line(0,1){20}}
% \put(110,80){\line(0,1){20}}
% \put(77,47){$1$}
% \put(97,87){$0$}
% \end{picture}
% \]
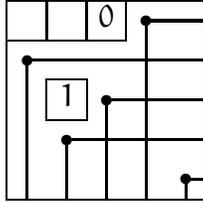
\[
\begin{picture}(150,75)
% \put(-50,30){$G(\pi),\ D(\pi)=$}
\put(37.5,0){\makebox[0pt][l]{\framebox(75,75)}}
\thicklines
\put(45,52.5){\circle*{4}}
\put(45,52.5){\line(1,0){67.5}}
\put(45,52.5){\line(0,-1){52.5}}

\put(60,22.5){\circle*{4}}
\put(60,22.5){\line(1,0){52.5}}
\put(60,22.5){\line(0,-1){22.5}}
\put(90,67.5){\circle*{4}}
\put(90,67.5){\line(1,0){22.5}}
\put(90,67.5){\line(0,-1){52.5}}
\put(75,37.5){\circle*{4}}
\put(75,37.5){\line(1,0){37.5}}
\put(75,37.5){\line(0,-1){37.5}}
\put(105,7.5){\circle*{4}}
\put(105,7.5){\line(1,0){7.5}}
\put(105,7.5){\line(0,-1){7.5}}
\thinlines
\put(52.5,30){\makebox[0pt][l]{\framebox(15,15)}}
\thicklines
\put(90,22.5){\line(1,0){7.5}}
\put(90,22.5){\line(0,-1){22.5}}
\thinlines
\put(37.5,60){\line(1,0){45}}
\put(67.5,60){\line(0,1){15}}
\put(52.5,60){\line(0,1){15}}
\put(82.5,60){\line(0,1){15}}
\put(57.25,35.25){$1$}
\put(72,65.25){$0$}
\end{picture}
\]
The diagram consists of two connected components; we also record the
value of the rank array~$r(\pi)$ on the essential set.

The matrix Schubert variety $\ol X_\pi$ is the set of $5\times 5$
matrices $Z$ such that $z_{11}=z_{12}=z_{13}=0$ and whose upper left
$3\times 2$ submatrix $Z_{3\times 2}$ has rank at most~1; all other
rank conditions on $Z \in \ol X_\pi$ follow from these.  Using only
the essential minors,
\begin{equation*}\label{eqn:ess_gen3142}
  I_\pi = \<z_{11}, \ z_{12}, \ z_{13}, \ z_{11}z_{22}-z_{21}z_{12}, \
  z_{11}z_{32}- z_{31}z_{12}, \ z_{21}z_{32}-z_{31}z_{22}\>.
\end{equation*}
These generators form a Gr\"obner basis under any diagonal term order,
i.e., one that picks $z_{11}$, $z_{12}$, $z_{13}$, $z_{11}z_{22}$,
$z_{11}z_{32}$, and $z_{21}z_{32}$ as the leading terms.  This
statement is an instance of the main result of this section,
Theorem~\ref{thm:grobnervex}.
\end{Example}

\subsection{Vexillary permutations}\label{ssec:vexillary}

Since we shall be interested in vexillary matrix Schubert varieties,
we collect in this subsection some results on vexillary permutations.

A permutation $\pi$ is called \dfn{vexillary} or \dfn{2143-avoiding}
if there do not exist integers $a<b<c<d$ such that $\pi(b) < \pi(a) <
\pi(d) < \pi(c)$.  Fulton characterized these permutations in terms of
their essential sets: no element $(p,q)$ of the essential set
is \dfn{strictly northwest} of another element $(i,j)$, meaning $p<i$ and $q<j$.
To each vexillary permutation
$\pi\in S_n$ we associate the partition $\lambda(\pi)$ whose parts are
the numbers of boxes in the rows of~$D(\pi)$, sorted into weakly
decreasing order.  For example, the permutation $\pi$
in Example~\ref{exa:3142} is vexillary, and $\lambda(\pi)=(3,1)$.

Every Grassmannian permutation (see the Introduction) is vexillary; in
fact, a permutation is Grassmannian if and only if its essential set
is contained in one row, necessarily the last nonempty row of the
diagram.  In this case, $\lambda(\pi)$ simply lists the number of
boxes in the rows of~$D(\pi)$, read bottom-up.

% Fulton has given a simple characterization of vexillary permutations
% in terms of the essential set \cite{fulton:92}.
In Corollary~\ref{cor:vex_equiv}, we will collect the characterizations
of vexillary permutations that we will need, one based on this.

\begin{Lemma} \label{lemma:pre_vex_equiv}
Let $\pi \in S_n$ and $(i,j) \in D(\pi)$.  There exists $(p,q) \in
\Ess(\pi)$ with $p < i$ and $q < j$ if and only if the $\rr ij\pi$
dots of $\pi_{(i-1) \times (j-1)}$ do not form a diagonal of
size~$\rr ij\pi$.
\end{Lemma}
\begin{proof}
This is straightforward from the definition of essential set.
\end{proof}

\begin{Corollary}[\cite{fulton:92}] \label{cor:vex_equiv}
The following are equivalent for a permutation $\pi\in S_n$.
\begin{alphlist}
\item
$\pi$ is vexillary.
\item
There do not exist $(p,q) \in \Ess(\pi)$ and $(i,j) \in \Ess(\pi)$
with $p < i$ and $q < j$.
\item
For all $(p,q)\in D(\pi)$, the $\rr pq\pi$ dots of~$\pi_{(p-1) \times
(q-1)}$ form a diagonal of size $\rr pq\pi$.
\end{alphlist}
\end{Corollary}
\begin{proof}
The equivalence of (a) and~(b) comes from \cite[Section~9]{fulton:92}.
Since every element of~$D(\pi)$ has an element of~$\Ess(\pi)$ to the
southeast of it in its connected component, the equivalence of (b)
and~(c) comes from Lemma~\ref{lemma:pre_vex_equiv}.%
\end{proof}

For any permutation~$\pi$, call a box $(p,q)\in D(\pi)$ \dfn{accessible} if 
$\rr pq\pi \neq 0$ and no boxes other than $(p,q)$ itself lie weakly to its
southeast in~$D(\pi)$.  In particular, $(p,q)\in \Ess(\pi)$.  Our next
goal is to define, for a vexillary~$\pi$ and an accessible box~$(p,q)$,
two new vexillary permutations $\pi_P$ and~$\pi_C$.

For any permutation $\pi$, each connected component of~$D(\pi)$ has
a unique northwest corner $(a,b)$. If $\pi$ is vexillary and $(a,b)\neq (1,1)$,
then there is a dot of~$\pi$ at $(a-1,b-1)$, because
Corollary~\ref{cor:vex_equiv} prevents pairs of dots of~$\pi$
weakly northwest of $(i,j)$ from forming antidiagonals.  
Let $(t,\pi(t))=(a-1,b-1)$
be the dot of~$\pi$ adjacent to the northwest corner of the connected
component of~$D(\pi)$ containing $(p,q)$.  Now set
\begin{equation}\label{eqn:pi_stuff}
  \pi_P=\pi \circ (p,\pi^{-1}(q)) \quad\mbox{and}\quad \pi_C=\pi_P
  \circ (t,p),
\end{equation}
where the composition $\mbox{}\circ(i,j)$ with the transposition
$(i,j)$ results in switching rows~$i$ and~$j$.  Denote the
corresponding rank matrices by $\rrr{P} = \rrr{\pi_P}$ and $\rrr{C} =
\rrr{\pi_C}$.

\begin{Example}\label{example:pi_CP}
The permutation $\pi=\left(\begin{array}{ccccccccc}
1 & 2 & 3 & 4 & 5 & 6 & 7 & 8 & 9 \\
8 & 7 & 1 & 6 & 2 & 9 & 5 & 3 & 4  
\end{array}\right)$
is vexillary.  Its dot-matrix and diagram are combined below:
% \[
% \begin{picture}(180,180)
% \put(0,0){\makebox[0pt][l]{\framebox(180,180)}}
% \put(150,170){\circle*{4}}
% \thicklines
% \put(150,170){\line(1,0){30}}
% \put(150,170){\line(0,-1){170}}
% \put(130,150){\circle*{4}}
% \put(130,150){\line(1,0){50}}
% \put(130,150){\line(0,-1){150}}
% \put(10,130){\circle*{4}}
% \put(10,130){\line(1,0){170}}
% \put(10,130){\line(0,-1){130}}
% 
% \put(110,110){\circle*{4}}
% \put(110,110){\line(1,0){70}}
% \put(110,110){\line(0,-1){110}}
% 
% \put(30,90){\circle*{4}}
% \put(30,90){\line(1,0){150}}
% \put(30,90){\line(0,-1){90}}
% 
% \put(170,70){\circle*{4}}
% \put(170,70){\line(1,0){10}}
% \put(170,70){\line(0,-1){70}}
% 
% \put(90,50){\circle*{4}}
% \put(90,50){\line(1,0){90}}
% \put(90,50){\line(0,-1){50}}
% 
% \put(50,30){\circle*{4}}
% \put(50,30){\line(1,0){130}}
% \put(50,30){\line(0,-1){30}}
% 
% \put(70,10){\circle*{4}}
% \put(70,10){\line(1,0){110}}
% \put(70,10){\line(0,-1){10}}
% 
% \thinlines
% 
% \put(0,140){\makebox[0pt][l]{\framebox(120,40)}}
% \put(0,160){\makebox[0pt][l]{\framebox(140,20)}}
% \put(20,140){\line(0,1){40}}
% \put(40,140){\line(0,1){40}}
% \put(60,140){\line(0,1){40}}
% \put(80,140){\line(0,1){40}}
% \put(100,140){\line(0,1){40}}
% 
% \put(20,100){\makebox[0pt][l]{\framebox(80,20)}}
% \put(40,100){\line(0,1){20}}
% \put(60,100){\line(0,1){20}}
% \put(80,100){\line(0,1){20}}
% 
% \put(40,40){\makebox[0pt][l]{\framebox(40,40)}}
% \put(40,60){\line(1,0){60}}
% \put(100,60){\line(0,1){20}}
% \put(80,80){\line(1,0){20}}
% \put(60,40){\line(0,1){40}}
% \end{picture}
% \]
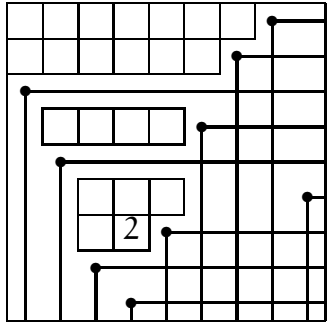
\[
\begin{picture}(120,120)
\put(0,0){\makebox[0pt][l]{\framebox(120,120)}}
\put(100,113.33){\circle*{4}}
\thicklines
\put(100,113.33){\line(1,0){20}}
\put(100,113.33){\line(0,-1){113.33}}
\put(86.67,100){\circle*{4}}
\put(86.67,100){\line(1,0){33.33}}
\put(86.67,100){\line(0,-1){100}}
\put(6.67,86.67){\circle*{4}}
\put(6.67,86.67){\line(1,0){113.33}}
\put(6.67,86.67){\line(0,-1){86.67}}

\put(73.33,73.33){\circle*{4}}
\put(73.33,73.33){\line(1,0){46.67}}
\put(73.33,73.33){\line(0,-1){73.33}}

\put(20,60){\circle*{4}}
\put(20,60){\line(1,0){100}}
\put(20,60){\line(0,-1){60}}

\put(113.33,46.67){\circle*{4}}
\put(113.33,46.67){\line(1,0){6.67}}
\put(113.33,46.67){\line(0,-1){46.67}}

\put(60,33.33){\circle*{4}}
\put(60,33.33){\line(1,0){60}}
\put(60,33.33){\line(0,-1){33.33}}

\put(33.33,20){\circle*{4}}
\put(33.33,20){\line(1,0){86.67}}
\put(33.33,20){\line(0,-1){20}}

\put(46.67,6.67){\circle*{4}}
\put(46.67,6.67){\line(1,0){73.33}}
\put(46.67,6.67){\line(0,-1){6.67}}

\thinlines

\put(0,106.67){\makebox[0pt][l]{\framebox(93.33,13.33)}}
\put(0,93.33){\makebox[0pt][l]{\framebox(80,26.67)}}
\put(13.33,93.33){\line(0,1){26.67}}
\put(26.67,93.33){\line(0,1){26.67}}
\put(40,93.33){\line(0,1){26.67}}
\put(53.33,93.33){\line(0,1){26.67}}
\put(66.67,93.33){\line(0,1){26.67}}

\put(13.33,66.67){\makebox[0pt][l]{\framebox(53.33,13.33)}}
\put(26.67,66.67){\line(0,1){13.33}}
\put(40,66.67){\line(0,1){13.33}}
\put(53.33,66.67){\line(0,1){13.33}}

\put(26.67,26.67){\makebox[0pt][l]{\framebox(26.67,26.67)}}
\put(26.67,40){\line(1,0){40}}
%\put(100,60){\line(0,1){20.25}}
\put(66.67,40){\line(0,1){13.5}}
%\put(80.15,80.2){\line(1,0){20}}
\put(53.43,53.58){\line(1,0){13.34}}
\put(40,26.67){\line(0,1){26.67}}
\put(43.75,30.75){$2$}
\end{picture}
\]
The box $(p,q)=(7,4)$, which is marked with $\rr 74\pi = 2$, is
accessible, and the dot
% given by Lemma~\ref{lemma:diag_adj}.
immediately northwest of its connected component is
$(t,\pi(t))=(5,2)$.  Therefore
\[\pi_P=\left(\begin{array}{ccccccccc}
1 & 2 & 3 & 4 & 5 & 6 & 7 & 8 & 9 \\
8 & 7 & 1 & 6 & 2 & 9 & 4 & 3 & 5  
\end{array}\right)
\mbox{ and \ \ }
\pi_C=\left(\begin{array}{ccccccccc}
1 & 2 & 3 & 4 & 5 & 6 & 7 & 8 & 9 \\
8 & 7 & 1 & 6 & 4 & 9 & 2 & 3 & 5  
\end{array}\right).
\]
These correspond respectively to 
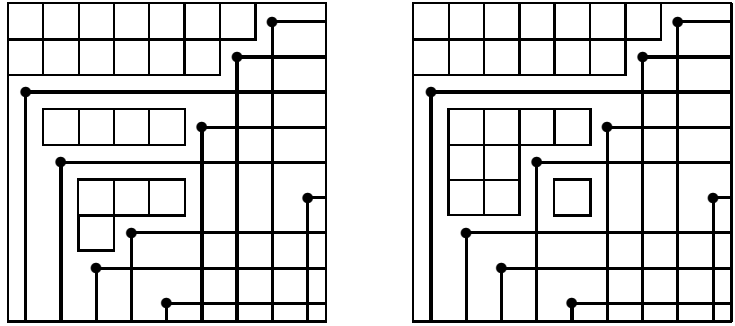
\[
\begin{picture}(273.33,120)
\put(153.33,0){\makebox[0pt][l]{\framebox(120,120)}}
\put(253.33,113.33){\circle*{4}}
\thicklines
\put(253.33,113.33){\line(1,0){20}}
\put(253.33,113.33){\line(0,-1){113.33}}
\put(240,100){\circle*{4}}
\put(240,100){\line(1,0){33.33}}
\put(240,100){\line(0,-1){100}}
\put(160,86.67){\circle*{4}}
\put(160,86.67){\line(1,0){113.33}}
\put(160,86.67){\line(0,-1){86.67}}

\put(226.67,73.33){\circle*{4}}
\put(226.67,73.33){\line(1,0){46.67}}
\put(226.67,73.33){\line(0,-1){73.33}}

\put(200,60){\circle*{4}}
\put(200,60){\line(1,0){73.33}}
\put(200,60){\line(0,-1){60}}

\put(266.67,46.67){\circle*{4}}
\put(266.67,46.67){\line(1,0){6.67}}
\put(266.67,46.67){\line(0,-1){46.67}}

\put(173.33,33.33){\circle*{4}}
\put(173.33,33.33){\line(1,0){100}}
\put(173.33,33.33){\line(0,-1){33.33}}

\put(186.67,20){\circle*{4}}
\put(186.67,20){\line(1,0){86.67}}
\put(186.67,20){\line(0,-1){20}}

\put(213.33,6.67){\circle*{4}}
\put(213.33,6.67){\line(1,0){60}}
\put(213.33,6.67){\line(0,-1){6.67}}

\thinlines

\put(153.33,93.33){\makebox[0pt][l]{\framebox(80,26.67)}}
\put(153.33,106.67){\makebox[0pt][l]{\framebox(93.33,13.33)}}
\put(166.67,93.33){\line(0,1){26.67}}
\put(180,93.33){\line(0,1){26.67}}
\put(193.33,93.33){\line(0,1){26.67}}
\put(206.67,93.33){\line(0,1){26.67}}
\put(220,93.33){\line(0,1){26.67}}

\put(166.67,66.67){\makebox[0pt][l]{\framebox(53.33,13.33)}}
\put(180,66.67){\line(0,1){13.33}}
\put(193.33,66.67){\line(0,1){13.33}}
\put(206.67,66.67){\line(0,1){13.33}}

\put(166.67,40){\line(1,0){26.67}}
\put(166.45,39.8){\line(0,1){26.67}}
\put(193.33,40){\line(0,1){26.67}}
\put(180,40){\line(0,1){26.67}}
\put(166.67,53.33){\line(1,0){26.67}}
\put(206.67,40){\makebox[0pt][l]{\framebox(13.4,13.33)}}

\put(0,0){\makebox[0pt][l]{\framebox(120,120)}}
\put(100,113.33){\circle*{4}}
\thicklines
\put(100,113.33){\line(1,0){20}}
\put(100,113.33){\line(0,-1){113.33}}
\put(86.67,100){\circle*{4}}
\put(86.67,100){\line(1,0){33.33}}
\put(86.67,100){\line(0,-1){100}}
\put(6.67,86.67){\circle*{4}}
\put(6.67,86.67){\line(1,0){113.33}}
\put(6.67,86.67){\line(0,-1){86.67}}

\put(73.33,73.33){\circle*{4}}
\put(73.33,73.33){\line(1,0){46.67}}
\put(73.33,73.33){\line(0,-1){73.33}}

\put(20,60){\circle*{4}}
\put(20,60){\line(1,0){100}}
\put(20,60){\line(0,-1){60}}

\put(113.33,46.67){\circle*{4}}
\put(113.33,46.67){\line(1,0){6.67}}
\put(113.33,46.67){\line(0,-1){46.67}}

\put(46.67,33.33){\circle*{4}}
\put(46.67,33.33){\line(1,0){73.33}}
\put(46.67,33.33){\line(0,-1){33.33}}

\put(33.33,20){\circle*{4}}
\put(33.33,20){\line(1,0){86.67}}
\put(33.33,20){\line(0,-1){20}}

\put(60,6.67){\circle*{4}}
\put(60,6.67){\line(1,0){60}}
\put(60,6.67){\line(0,-1){6.67}}

\thinlines

\put(0,93.33){\makebox[0pt][l]{\framebox(80,26.67)}}
\put(0,106.67){\makebox[0pt][l]{\framebox(93.33,13.33)}}
\put(13.33,93.33){\line(0,1){26.67}}
\put(26.67,93.33){\line(0,1){26.67}}
\put(40,93.33){\line(0,1){26.67}}
\put(53.33,93.33){\line(0,1){26.67}}
\put(66.67,93.33){\line(0,1){26.67}}

\put(13.33,66.67){\makebox[0pt][l]{\framebox(53.33,13.33)}}
\put(26.67,66.67){\line(0,1){13.33}}
\put(40,66.67){\line(0,1){13.33}}
\put(53.33,66.67){\line(0,1){13.33}}

\put(26.67,40){\makebox[0pt][l]{\framebox(40,13.33)}}
\put(40,26.67){\line(0,1){26.67}}
\put(26.67,26.67){\line(1,0){13.33}}
\put(26.5,26.45){\line(0,1){13.33}}
\put(53.33,40){\line(0,1){13.33}}
\end{picture}
\]
\end{Example}

\begin{Lemma}\label{lem:CPmatrices}
Let $\pi$ be a vexillary permutation, $\rrr\pi$ its rank matrix, and
$(p,q)$ an accessible box for $\pi$.  Then the following hold.
\begin{alphlist}
\item
$D(\pi_P)=D(\pi)\minus \{(p,q)\}$.
\item
$D(\pi_C)$ is obtained by moving diagonally northwest by one step the
rectangle consisting of boxes of $D(\pi)$ weakly northwest of $(p,q)$
and in its connected component.\vspace{1ex}
\item
$\rr ijP =
\left\{\begin{array}{ll}
  \rr ij\pi + 1 & \hbox{if } p\leq i\leq \pi^{-1}(q)-1 \hbox{ and }
  q\leq j\leq \pi(p)-1
\\
  \rr ij\pi & \hbox{otherwise.}
\end{array}\right.$\vspace{2ex}
\item
$\rr ijC = 
\left\{\begin{array}{ll}
  \rr ij\pi - 1 &\hbox{if } t\leq i\leq p-1 \hbox{ and }\pi(t)\leq
  j\leq q-1
\\
  \rr ij\pi + 1 &\hbox{if } p\leq i\leq \pi^{-1}(q)-1 \hbox{ and }
  q\leq j\leq \pi(p)-1
\\
  \rr ij\pi & \mbox{otherwise.}
\end{array}\right.$\vspace{1ex}
\item
$\Ess(\pi) \minus \{(p,q)\} \subseteq \Ess(\pi_P)$ and
$\Ess(\pi_P)\minus \Ess(\pi)\subseteq \{(p-1,q),(p,q-1)\}$.
\item
$\Ess(\pi_C) = (\Ess(\pi)\minus\{(p,q)\}) \cup \{(p-1,q-1)\}$.
\item
$\pi_P$ and $\pi_C$ are vexillary permutations.
\end{alphlist}
\end{Lemma}
\begin{proof}
Parts~(a)--(f) are straightforward to check from the definitions.
Part~(g) follows easily from parts (a) and~(b) combined with the
equivalence of Corollary~\ref{cor:vex_equiv}(a) and
Corollary~\ref{cor:vex_equiv}(b).
\end{proof}

\def\5{ .}
\def\6{ .}
\def\7{ .}
\def\8{ .}
\def\9{ .}

\begin{Example}
Continuing Example~\ref{example:pi_CP}, omitting the $2$'s to better
see the shapes, and omitting the $5$'s, $6$'s, $7$'s, $8$'s, and $9$'s
(which don't change), we have
\[
\rrr\pi=\left(
\begin{matrix}
 0 &  0 &  0 &  0 &  0 &  0 &  0 &  1 &  1 \\
 0 &  0 &  0 &  0 &  0 &  0 &  1 &  . &  .\\
 1 &  1 &  1 &  1 &  1 &  1 &  . &  3 &  3 \\
 1 &  1 &  1 &  1 &  1 &  . &  3 &  4 &  4 \\
 1 &  . &  . &  . &  . &  3 &  4 & \5 & \5 \\
 1 &  . &  . &  . &  . &  3 &  4 & \5 & \6 \\
 1 &  . &  . &  . &  3 &  4 & \5 & \6 & \7 \\
 1 &  . &  3 &  3 &  4 & \5 & \6 & \7 & \8 \\
 1 &  . &  3 &  4 & \5 & \6 & \7 & \8 & \9 
\end{matrix}
\right),
\]
\[\rrr{P}=\left(
\begin{matrix}
 0 &  0 &  0 &  0 &  0 &  0 &  0 &  1 &  1 \\
 0 &  0 &  0 &  0 &  0 &  0 &  1 &  . &  .\\
 1 &  1 &  1 &  1 &  1 &  1 &  . &  3 &  3 \\
 1 &  1 &  1 &  1 &  1 &  . &  3 &  4 &  4 \\
 1 &  . &  . &  . &  . &  3 &  4 & \5 & \5 \\
 1 &  . &  . &  . &  . &  3 &  4 & \5 & \6 \\
 1 &  . &  . &  3 &  3 &  4 & \5 & \6 & \7 \\
 1 &  . &  3 &  4 &  4 & \5 & \6 & \7 & \8 \\
 1 &  . &  3 &  4 & \5 & \6 & \7 & \8 & \9 
\end{matrix}
\right), \quad\hbox{and}\quad
\rrr{C}=\left(
\begin{matrix}
 0 &  0 &  0 &  0 &  0 &  0 &  0 &  1 &  1 \\
 0 &  0 &  0 &  0 &  0 &  0 &  1 &  . &  .\\
 1 &  1 &  1 &  1 &  1 &  1 &  . &  3 &  3 \\
 1 &  1 &  1 &  1 &  1 &  . &  3 &  4 &  4 \\
 1 &  1 &  1 &  . &  . &  3 &  4 & \5 & \5 \\
 1 &  1 &  1 &  . &  . &  3 &  4 & \5 & \6 \\
 1 &  . &  . &  3 &  3 &  4 & \5 & \6 & \7 \\
 1 &  . &  3 &  4 &  4 & \5 & \6 & \7 & \8 \\
 1 &  . &  3 &  4 & \5 & \6 & \7 & \8 & \9 
\end{matrix}
\right).
\]
These rank arrays $\rrr{P}$ and~$\rrr{C}$ are in agreement with
Lemma~\ref{lem:CPmatrices}.%
\end{Example}

Finally, we see how to get vexillary permutations from Grassmannian
ones.

\begin{Lemma} \label{lemma:chain}
Let $\pi \in S_n$ be vexillary permutation with largest descent
position~$k$.  Then there exists some $N \geq n$ and a sequence of
vexillary permutations
\[
  \sigma_1,\sigma_2,\ldots,\sigma_t = \pi
\]
in~$S_N$, all with largest descent position~$k$, such that $\sigma_1$
is \pagebreak[2] Grassmannian and $\sigma_{i+1}=(\sigma_i)_C$ for
$1\leq i\leq t-1$.  The permutation $\tpi = \sigma_1$ is uniquely
determined~by~$\pi$; in fact,
\[
  \Ess(\tpi) = \{(k,k-p+q) \mid (p,q) \in \Ess(\pi)\}
  \quad\text{and}\quad \rr{k,}{k-p+q}\tpi = k-p+\rr pq\pi.
\]
\end{Lemma}
\begin{proof}
There is nothing to prove if $\pi$ is Grassmannian.  Otherwise we
construct a vexillary permutation~$\sigma$ as follows.  Find the
second largest descent $i < k$ of~$\pi$.  Hence the rightmost
box~\mbox{$(i,j)$} of~$D(\pi)$ in row~$i$ lies in~$\Ess(\pi)$.  Since
$\pi$ is assumed to be vexillary, all boxes of~$\Ess(\pi)$ to the
north of~\mbox{$(i,j)$} are weakly to the east of column~$j$.  Find
the northmost box in~$\Ess(\pi)$ that is in column~$j$,
say~\mbox{$(h,j)$}.

Since $(h,j) \in \Ess(\pi)$, there are dots $(h+1,q)$ and $(p,j+1)$
of~$\pi$ satisfying $q \leq j$ and $p \leq h$.  Since $\pi$ is
vexillary, there are no boxes of the diagram~$D(\pi)$ strictly
southeast of~\mbox{$(h,j)$}. This implies that there is a unique
northwestmost dot $(c,\pi(c))$ of $\pi$ strictly southeast
of~\mbox{$(h,j)$}; that is, no two dots of~$\pi$ strictly southeast
of~\mbox{$(h,j)$} form an antidiagonal.

Now let $\sigma = \pi \circ (p,h+1) \circ (h+1,c)$.  Then one checks
that $D(\sigma)$ is obtained from~$D(\pi)$ by moving the boxes in the
rectangle in row~$p$ to~$h$ and columns~$q$ and~$j$ southeast by one
unit.  Our choices guarantee that $\sigma$ is still vexillary, by
looking at $\Ess(\sigma)$ and using Corollary~\ref{cor:vex_equiv}.
Moreover, $\lambda(\sigma) = \lambda(\pi)$ because the number of rows
with any given number of boxes of the diagram is the same for both
$\pi$ and~$\sigma$.  Also, $\sigma$ still has its last descent at
position~$k$.  Most importantly, $(h+1,j+1)\in \Ess(\sigma)$ is
accessible, and $\pi = \sigma_C$.
          
Note in particular that the boxes of~$\sigma$ are further south than
those of~$\pi$ (and strictly so in at least one case).  Therefore,
repeated application of the above construction successively moves the
boxes of $D(\pi)$ south, but always above row~$k$.  So this gives a
chain of permutations starting with $\pi$ and eventually ending with a
Grassmannian permutation in some $S_N$ with $N \geq n$.

The final sentence of the lemma follows from parts~(d) and~(f) of
Lemma~\ref{lem:CPmatrices}.
\end{proof}

\subsection{Diagonal Gr\"obner bases}\label{ssec:diagGrobner}

In \cite{grobGeom} it was proved that for any permutation~$\pi$, the
essential minors form a Gr\"obner basis of the Schubert determinantal
ideal~$I_\pi$ under any \emph{anti}\/diagonal term order.  We prove
here a complementary result for diagonal term orders, assuming that
the permutation~$\pi$ is vexillary (Theorem~\ref{thm:grobnervex}); in
Section~\ref{section:sharp} we explain the sense in which this
diagonal Gr\"obner basis result is sharp.  Our proof will be an
application of Theorem~\ref{thm:algGVD}.
% \comment{EM: The following has already been said in the
% Introduction; I don't see a need to repeat it here:}
% The class of vexillary Schubert determinantal ideals is known to
% coincide with the class of (one-sided) ladder determinantal ideals;
% see, for example, \cite[Section~2.4]{grobGeom} and the references
% therein.

\begin{Theorem}\label{thm:grobnervex}
If $\pi \in S_n$ is a vexillary permutation, then the essential minors
of~$I_\pi$ constitute a Gr\"obner basis with respect to any diagonal
term order.
\end{Theorem}

The proof will be by an induction based on Lemma~\ref{lemma:chain} and
the following, which connects the notation of
Theorem~\ref{thm:algGVD} with Eq.~(\ref{eqn:pi_stuff}).
Geometrically, Proposition~\ref{p:inductionstep} concerns the
decomposition $M_n = H\oplus L$ where $H$ consists of all matrices
satisfying % the rank conditions imposed by a vexillary
% permutation~$\pi$ together with the condition that 
$z_{pq}=0$ for some
fixed accessible box $(p,q)\in \Ess(\pi)$, and $L$ is the $1$-dimensional
space of matrices with all $z_{ij}$ vanishing except $z_{pq}$.

\begin{Proposition}\label{p:inductionstep}
  Let $\pi$ be a vexillary permutation and $(p,q)\in \Ess(\pi)$ an
  accessible box.  Fix a diagonal term order~$\prec$ on~$\kk[\zz]$.
  Suppose that the essential minors generating $I=I_\pi$ form a
  Gr\"obner basis with respect to~$\prec$,
% (This will be seen later to be automatic.)
% Let $(p,q)$ be an accessible box and let $M_{n}=H\oplus L$
% where $H$ are those matrices whose $(p,q)$ entry is zero.
  set $z_{pq} = y$, and let $I'$, $C$, $P$, and their respective
  generators be as in Theorem~\ref{thm:algGVD}.  Then the
  following hold.
  \begin{alphlist}
  \item
  $C = I_{\pi_C}$ and $P = I_{\pi_P} + \< z_{pq}\>$, where $\pi_C$ and
  $\pi_P$ are defined by Eq.~(\ref{eqn:pi_stuff}).
  \item
  The generators of $C$ and $P$ are Gr\"obner bases for these ideals
  under~$\prec$.
  \item
  The essential minors generating $I_{\pi_C}$ form a Gr\"obner
  basis under~$\prec$.
  \item
  $I' = C \cap P$ is a geometric vertex decomposition.
  \end{alphlist}
\end{Proposition}
\begin{proof}
Separate the set of essential minors generating $I$ into sets~$V$
and~$U$ according to whether they do or do not involve~$z_{pq}$.
Since $(p,q)$ is accessible, the array of $\zz$-variables appearing in
an essential minor that involves~$z_{pq}$ has southeast
corner~$z_{pq}$.
  
For each minor in~$V$, put into a new set $V'$ the one-smaller minor
obtained by removing the last row and column.  Define $z_{pq} V'$ by
multiplying all elements of~$V'$ by~$z_{pq}$.  Then, by definition
of Gr\"obner basis,
$I'$ is generated by $U \cup (z_{pq} V')$, whereas $C$ is generated by
$U \cup V'$ and $P$ is generated by $U \cup \{z_{pq}\}$.

It follows from parts~(c) and~(e) of Lemma~\ref{lem:CPmatrices} that
the essential minors for $\pi_P$ are precisely the determinants
in~$U$.  Hence the statement about~$P$ in part~(a) follows.

Next we check that $C\subseteq I_{\pi_C}$.  Each minor in~$U$ comes
from a (possibly inessential) rank condition from $\rrr\pi$
corresponding to either
\begin{romanlist}
\item
a box of $D(\pi)\minus \{(p,q)\}$ in the same connected component as
$(p,q)$ and in row~$p$ or column~$q$, or
\item
a box in $\Ess(\pi)\minus\{(p,q)\}$.
\end{romanlist}
In either case, Lemma~\ref{lem:CPmatrices}(d) shows that $\rrr{\pi_C}$
and~$\rrr\pi$ are equal at these positions.  Hence each minor in~$U$
lies in~$I_{\pi_C}$.  Also, any minor in~$V$ arises from the rank
condition corresponding to $(p,q) \in \Ess(\pi)$, and the associated
minor in~$V'$ is an essential minor corresponding to $\rr
{p-1,}{q-1}C$, by parts~(d) and~(f) of Lemma~\ref{lem:CPmatrices}.
Thus we obtain $V'\subseteq I_{\pi_C}$, and so $U \cup V' \subseteq
I_{\pi_C}$.  Therefore $C \subseteq I_{\pi_C}$.

To show the other inclusion, one checks from
Lemma~\ref{lem:CPmatrices} that if an essential minor comes from
$(p-1,q-1)\in \Ess(\pi_C)$ then it lies in~$V'$; otherwise, it lies
in~$U$.  This concludes the proof of part~(a).

Theorem~\ref{thm:algGVD} guarantees that the generating sets
for~$C$ and~$P$ are again Gr\"obner bases.  Since minors have only
squarefree terms, the power of~$z_{pq}$ in each generator is at
most~$1$, so the theorem also tells us we have a geometric vertex
decomposition.  This proves parts~(b) and~(d).

It remains to check part~(c).  Although $I_{\pi_C} = C$, there are
minors in $U \cup V'$ that are not essential minors for $I_{\pi_C}$.
These are exactly those minors in~$U$ arising from~(i) above.  Thus we
wish to show that we can remove these inessential minors and still
have a Gr\"obner basis.  Since~$I_{\pi_C}$ is generated by its
essential minors, it suffices to check that the leading term of each
minor from~(i) is divisible by the leading term of an essential minor
in~$I_{\pi_C}$.  For this, note that removing the last row and column
from one of these inessential minors yields a minor of smaller degree
that is an essential minor of~$I_{\pi_C}$ arising
from~$\rr{p-1,}{q-1}C$.%
\end{proof}

\begin{proofof}{Theorem~\ref{thm:grobnervex}}
% The property of being a Gr\"obner basis does not change under the
% presence of extra variables that do not appear in any of the
% polynomials.
Construct a chain of vexillary permutations as in
Lemma~\ref{lemma:chain}.  If $\pi = \sigma_1$, so $\pi$ is Grassmannian, then
the result is proved in \cite{CGG,Stu90}.  This case can also be
derived from \cite[Theorem~B]{grobGeom}, which proves the Gr\"obner
basis property for antidiagonal term orders, because the set of
essential minors for a Grassmannian permutation with unique
descent~$k$ is invariant under the permutation of the variables
induced by reversing the top~$k$ rows of~$\zz$ (since the essential
set of~$\pi$ lies entirely in the $k^\mathrm{th}$ row of $D(\pi)$).
If $\pi \neq \sigma_1$ then the desired statement holds by induction
on the length of the sequence in Lemma~\ref{lemma:chain}, using
Proposition~\ref{p:inductionstep}(c).%
\end{proofof}\vspace{-1ex}

%%%%%%%%%%%%%%%%%%%%%%%%%%%%%%%%%%%%%%%%%%%%%%%%%%%%%%%%%%%%%%%%%%%%%%
\section{Stanley--Reisner ideals and subword complexes}\label{section:subword}
%%%%%%%%%%%%%%%%%%%%%%%%%%%%%%%%%%%%%%%%%%%%%%%%%%%%%%%%%%%%%%%%%%%%%%

%\comment{Define here the subword complex $\Gamma_\pi$ 
%  relevant to a vexillary permutation, so called because in [KM]
%  they were $L_\pi$ and $\Gamma$ looks like $L$ upside down.}
Theorem~\ref{thm:grobnervex} shows that for a vexillary permutation
$\pi\in S_n$, the initial ideal $\init\,I_\pi$ is a squarefree
monomial ideal; this makes it the \dfn{Stanley--Reisner ideal} of a
simplicial complex whose vertices are $[n]^{2}=\{(p,q)\mid 1\leq
p,q\leq n\}$ and whose faces consist of those subsets~$F$ such that no
monomial from $\init\,I_\pi$ has support~$F$.  

In this section we prove this complex to be a subword complex~\cite{subword}; 
besides its intrinsic interest, we apply this fact toward
combinatorial formulae in the next section. Without introducing
subword complexes, though, our geometric technology is already
sufficient to prove that this complex is shellable.

\begin{Theorem}\label{thm:gvdshelling}
  Let $\pi$ be a vexillary permutation, and $\Gamma$ the simplicial complex
  whose Stanley--Reisner ideal is $\init\, I_\pi$. Then $\Gamma$ is shellable.
\end{Theorem}

\begin{proof}
  If $J$ is a Stanley--Reisner ideal in $\kk[\zz]$, let $\Gamma(J)$
  denote the corresponding simplicial complex on the vertex
  set~$[n]^2$, so $\Gamma = \Gamma(\init\, I_\pi)$.

  By Proposition~\ref{p:inductionstep} part (d), $\ol X_\pi$ has a
  geometric vertex decomposition into $\ol X_{\pi_P} \cap \{z_{pq} = 0\}$
  and $\ol X_{\pi_C}$, where $(p,q)$ is an accessible box of $\pi$.
  By Proposition~\ref{prop:inductivegvd}, 
  this geometric vertex decomposition of $\ol X_\pi$ degenerates to
  a geometric vertex decomposition of $\init\, I_\pi$ into
  $\init\, I_{\pi_P} + \<z_{pq}\>$ and $\init\, I_{\pi_C}$.
  As explained in Example~\ref{ex:SR}, this 
  gives an ordinary vertex decomposition of $\Gamma(\init\, I_\pi)$
  at the vertex $(p,q)$
  into the deletion $\Gamma(\init\, I_{\pi_P} + \<z_{pq}\>)$ and the
  cone $\Gamma(\init\, I_{\pi_C})$ on the link.
  
  Since the ideal $\init\, I_{\pi_P}$ doesn't involve the generator
  $z_{pq}$, the complex $\Gamma(\init\, I_{\pi_P})$ is the cone on
  $\Gamma(\init\, I_{\pi_P} + \<z_{pq}\>)$ at the vertex $(p,q)$. 
  In particular, one is shellable if and only if the other is.

  Using induction on
% e.g.\
  the position of the most southeastern box in the diagram,
  $\Gamma(\init\, I_{\pi_C})$ and $\Gamma(\init\, I_{\pi_P})$ are both
  shellable; hence $\Gamma(\init\, I_{\pi_P} + \<z_{pq}\>)$ is also.
  As noted in \cite{BP}, one can concatenate a shelling of the
  deletion $\Gamma(\init\, I_{\pi_P} + \<z_{pq}\>)$ and a shelling of
  the cone $\Gamma(\init\, I_{\pi_C})$ on the link to make a shelling
  of $\Gamma(\init\, I_\pi)$.
\end{proof}

We review some definitions about subword complexes
from~\cite{subword}; see also \cite[Section~16.5]{cca}, which covers
the generality here.  A \dfn{word} of size $t$ is an ordered sequence
$Q=(s_{i_1},\ldots,s_{i_t})$ of simple reflections $s_i=(i,i+1)\in
S_{m}$.  An ordered subsequence $P$ of~$Q$ is a \dfn{subword} of~$Q$.
In addition, $P$ \dfn{represents} $\rho\in S_m$ if the ordered product
of the simple reflections in~$P$ is a reduced decomposition of $\rho$.
We say that $P$ \dfn{contains} $\rho \in S_m$ if some subsequence
of~$P$ represents~$\rho$.  Let $\Delta(Q,\rho)$ denote the
\dfn{subword complex} with vertex set the simple reflections in~$Q$
and with faces given by the subwords $Q \minus P$ whose complements
$P$ \mbox{contain~$\rho$}.
% this \rho isn't the \tpi from Lemma:chain, is it?

% Using the above setup, we now define the subword complex $L_\pi$.
Using the above setup, we now define the subword complex $\Gamma_\pi$
for a given vexillary permutation $\pi \in S_n$.  Let $\mu(\pi) =
(\mu_1 \geq \mu_2 \geq \ldots \geq \mu_k > 0)$ be the partition with
the smallest Ferrers shape (in English notation, i.e. the largest part
along the top row) containing all of the boxes of~$D(\pi)$;
%i.e., for each row $j$ such that
%there are boxes of $D(\pi)$ weakly to the south, $\mu_j$ is 
%the index of the rightmost column occupied by any such box of $D(\pi)$. 
thus the shape of~$\mu(\pi)$
% consists of all locations that are weakly northwest of some box
% in~$\Ess(\pi)$.
is the union over $(p,q) \in \Ess(\pi)$ of the northwest $p \times q$
rectangles in the $n \times n$ grid.  Fill each box $(p,q)$ in the
shape of~$\mu(\pi)$ with the reflection~$s_{k-p+q}$.  Reading each row
of~$\mu(\pi)$ from \emph{right to left}, starting with the bottom row
and ending with the top row yields the word
\[
  Q = (
  \underbrace{s_{\mu_k},\ldots,s_1}_{\text{row }k},
  \underbrace{s_{1+\mu_{k-1}},\ldots,s_2}_{\text{row }k - 1},
  \ldots,
  \underbrace{s_{i+\mu_i},\ldots,s_{i+1}}_{\text{row } k-i},
  \ldots,
  \underbrace{s_{k-1+\mu_1},\ldots,s_k}_{\text{top row}}
  ).
\]
% There is a bijective correspondence between reflections in $Q$ and
% the shape of the ``reverse partition'' ${\hat
% \mu}(\pi)=(\mu_k,\ldots,\mu_1)$ given by filling the
% $i^\mathrm{th}$ row of this shape by $i,i+1,\ldots,\mu_i$ from left to
% right, top down, i.e., $Q$ is obtained by reading the rows of this
% filling right to left and top down.
% 
% Now set $N = \mu_1 + k$, and denote by $\tpi \in S_N$ the unique
% Grassmannian permutation of shape~$\lambda(\pi)$ with descent at
% position~$k$.
Let $\tpi$ be the Grassmannian permutation from
Lemma~\ref{lemma:chain}, and set $\Gamma_\pi = \Delta(Q,\tpi)$.

\begin{Example} \label{ex:Gamma}
The vexillary permutation $\pi=\left(\begin{array}{ccccc}
1 & 2 & 3 & 4 & 5 \\
4 & 1 & 3 & 2 & 5 
\end{array}\right)\in S_{5}$ from Example~\ref{exa:3142}
has $\mu(\pi)=(3,2,2)$ and $\tpi=
\left(\begin{array}{cccccc}
1 & 2 & 3 & 4 & 5 & 6\\
1 & 3 & 6 & 2 & 4 & 5 
\end{array}\right)\in S_N$ for $N = 6$.  The ambient word
for~$\Gamma_\pi$ is $Q=(s_2,s_1,s_3,s_2,s_5,s_4,s_3)$.
\end{Example}

\begin{Lemma} \label{l:mu}
Identify~$\mu(\pi)$ with its Ferrers shape.  Setting to $1 \in \kk$
all variables $z_{pq}$ for $(p,q)$ outside of~$\mu(\pi)$ takes any
generating set for\/ $\IN\,I_\tpi$ to a generating set for\/
$\IN\,I_\pi$.
\end{Lemma}
\begin{proof}
It suffices to check that (i)~every diagonal generator of~$\IN\,I_\pi$
can be obtained from some diagonal in~$\IN\,I_\tpi$ by setting the
variables outside of~$\mu(\pi)$ to~$1$, and (ii)~setting the given
variables to~$1$ in any diagonal from~$\IN\,I_\tpi$ yields an element
inside~$\IN\,I_\pi$.
% $\IN\,I_\tpi = \<\text{all diagonals of size } k-p+1+\rr pq\pi
% \text{ in } \zz_{k\times(k-p+q)}\text{ for all}(p,q)\in\Ess(\pi)\>$
% whereas
% $\init\,I_\pi = \<\text{all diagonals of size } 1+\rr pq\pi \text{
% in} \zz_{p \times q}\text{ for all }(p,q)\in\Ess(\pi)\>$.
For~(i), if $\delta$ is a diagonal generator of~$I_\pi$, say of degree
$1+\rr pq\pi$ using the variables $\zz_{p \times q}$ for $(p,q) \in
\Ess(\pi)$, then simply multiply~$\delta$ by the diagonal $z_{1+p,1+q}
\cdots z_{k,k-p+q}$ to get a diagonal generator of~$\IN\,I_\tpi$.
For~(ii), observe that any diagonal of size at least $j+1+\rr pq\pi$
using the variables $\zz_{(j+p)\times(j+q)}$ has at least $1+\rr
pq\pi$ of its variables weakly northwest of~$z_{pq}$, and take $j =
k-p$ for $(p,q) \in \Ess(\pi)$.%
\end{proof}

\begin{Theorem}\label{thm:grobnervex'}
Let $\pi\in S_n$ be a vexillary permutation.  With respect to any
diagonal term order, the initial ideal of~$I_\pi$ is the
Stanley--Reisner ideal for the subword complex~$\Gamma_\pi$.
\end{Theorem}
\begin{proof}
Let $\tpi$ be the Grassmannian permutation defined in Lemma~\ref{lemma:chain}.
If $\pi = \tpi$ then the result follows from \cite[Theorem~B and
Example~1.8.3]{grobGeom}, which proves the corresponding result for
the \emph{anti}\/diagonal initial ideal of~$I_\tpi$, because the set
of essential minors for $I_{\tpi}$ is invariant under the permutation of
the variables induced by reversing the top~$k$ rows of~$\zz$.

For general vexillary~$\pi$, let us write $\Gamma_\tpi = \Delta(\tilde
Q,\tpi)$, to distinguish the word $\tilde Q$ from~$Q$.  Given a
subword~$P$ of~$\tilde Q$, denote by~$\zz_P$ the set of variables
corresponding to the locations of the boxes occupied by~$P$.  The
previous paragraph yields the minimal prime decomposition $\IN\,I_\tpi
= \bigcap_P \<\zz_P\>$, the intersection being over subwords $P$
of~$\tilde Q$ representing~$\tpi$.  Since $\IN\,I_\tpi$ is a monomial
ideal, setting the variables outside of~$\mu(\pi)$ to $1 \in \kk$ and
omitting those intersectands $\<\zz_P\>$ that become the unit ideal
yields a prime decomposition, and it is a decomposition
of~$\IN\,I_\pi$ by Lemma~\ref{l:mu}.  (The intersections for
$\IN\,I_\tpi$ and $\IN\,I_\pi$ are taken in different polynomial
rings, but this is irrelevant here.)  The intersectand $\<\zz_P\>$ for
a subword $P \subseteq \tilde Q$ representing~$\tpi$ survives the
process of setting the variables outside of~$\mu(\pi)$ to~$1$ if and
only if~$P$ is actually a subword of~$Q$.%
% Therefore $\IN\,I_\pi = \bigcap_P \<\zz_P\>$, the intersection being
% taken over all subwords $P$ of~$\tilde Q$ that represent~$\tpi$ and
% happen also to be subwords of~$Q$.%
\end{proof}

In the next section, we will use a pictorial description
of~$\Gamma_\pi$, developed in \cite{FKyangBax,BB,grobGeom} (see
\cite[Chapter~16]{cca} for an exposition).  Let $\pi$ be vexillary,
and again consider the Grassmannian permutation $\tpi \in S_N$ with
descent at~$k$.  A~tiling of the
% Ferrers shape $\mu(\pi)$ by
$k \times N$ rectangle by \dfn{crosses} $\textcross$ and
\dfn{elbows}~$\textelbow$ is called a \dfn{pipe dream}.  (Warning: the
elbow tile here is mirror-reflected top-to-bottom from the above
references. This is traceable to our use of diagonal, rather than antidiagonal,
term orders.) It is often convenient to identify a pipe dream with its
set of $\textcross$ tiles, and to identify each subword~$P$ of~$Q$
with the pipe dream whose crosses lie at the positions occupied
by~$P$.  When $P$ represents~$\tpi$, so the number of crosses equals
the length of~$\tpi$, then $P$ is a \dfn{reduced} pipe dream (or an
``rc-graph'') for~$\tpi$.
% When drawing pipe dreams for $\tpi$ we omit drawing the ``sea of
% waves'' that appear in the lower triangular part of the diagram.
In this language, the initial ideal $\IN\,I_\pi$ is described as
follows, by \cite[Example~1.8.3]{grobGeom}.

\begin{Corollary}
$Q \minus P$ is a face of $\Gamma_\pi$ if and only if (i)~the crosses
in the pipe dream~$P$ lie in the Ferrers shape~$\mu(\pi)$, and
(ii)~$P$ contains some reduced pipe dream representing~$\tpi$.
\end{Corollary}

\begin{Example} \label{ex:P}
Recall the situation from Example~\ref{ex:Gamma}.  The following pipe
dream~$P$ corresponds to the subword $(s_2,\ \cdot\ ,\ \cdot\ ,\
\cdot\ ,s_5,s_4,s_3)$ of~$Q$, and $Q \minus P$ is a facet
of~$\Gamma_\pi$.
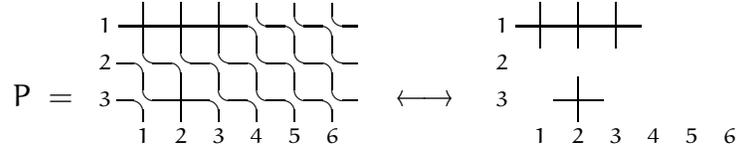
\[
\begin{array}{c}
  \\\\P\ =\\\\
\end{array}
\begin{array}{cccccccccccc}
  \petit1 &   \+   &   \+   &   \+   &  \rj   &   \rj  &  \rj   \\
  \petit2 &   \rj  &   \rj  &   \rj  &  \rj   &   \rj  &  \rj   \\
  \petit3 &   \rj  &   \+   &   \rj  &  \rj   &   \rj  &  \rj   \\
          &\perm1{}&\perm2{}&\perm3{}&\perm4{}&\perm5{}&\perm6{}\\
\end{array}
\begin{array}{c}
  \\\\\ \ \ \ \longleftrightarrow\ \ \\\\
\end{array}
\begin{array}{ccccccc}
  \petit1 & \+  & \+  & \+  &     &     &     \\
  \petit2 &     &     &     &     &     &     \\
  \petit3 &     & \+  &     &     &     &     \\
          &\perm1{}&\perm2{}&\perm3{}&\perm4{}&\perm5{}&\perm6{}\\
\end{array}
\]
By Example~\ref{exa:3142}, $\IN\,I_\pi = \<z_{11}, z_{12}, z_{13},
z_{21}z_{32}\> = \<z_{11}, z_{12}, z_{13}, z_{21}\> \cap \<z_{11},
z_{12}, z_{13}, z_{32}\>$.  Thus there is only one remaining facet $Q
\minus P'$ of~$\Gamma_\pi$, corresponding to the (reduced) pipe
dream~$P'$ obtained by moving the $\textcross$ at~$(3,2)$ diagonally
northwest to~$(2,1)$.

If we follow the three pipes from the left side to the bottom, they come
out in positions $1,3,6$. These are the positions of the ``up'' moves
in the word ``up, right, up, right, right, up'' describing a walk around
the partition $(3,1,0)$ from the southwest corner to the northeast.
\end{Example}

\begin{excise}{%
  The faces of $L_\pi$ can be represented in terms of \dfn{pipe
  dreams}, see, e.g., \cite{grobGeom}.  These are tilings of boxes in
  a shape by a \emph{crosses} $\textcross$ or an \emph{elbow joints}
  $\textelbow$; although we will often find it convenient to describe
  them only in terms of the crosses that they contain.  To each face
  $Q\minus P$ in $L_\pi$ we associate a pipe dream $\PP$ by placing
  crosses at the positions of the reverse partition $\mu(\pi)$ that
  correspond to $P$, placing this arrangement of crosses in the upper
  left corner of the $m\times m$ square and then filling the remaining
  boxes by elbows. This is a \dfn{pipe dream for $\tpi$}, see, e.g.,
  \cite{grobGeom}. When the number of crosses equals the length of
  $\tpi$, i.e., when $P$ represents $\tpi$, then $\PP$ is
  \dfn{reduced} or an ``RC-graph'', see \cite{FKyangBax,BB}. When
  drawing pipe dreams for $\tpi$ we omit drawing the ``sea of waves''
  that appear in the lower triangular part of the diagram.
  
  \begin{Example} \label{exa:base_pipe}
  The subword $P=(\cdot,\cdot,\cdot,s_2,s_5,s_4,s_3)$ corresponds to
  the pipe dream
  % $\PP$
  \[
  \PP = 
    \begin{array}{cccccccccccc}
            &\perm1{}&\perm2{}&\perm3{}&\perm4{}&\perm5{}&\perm6{}\\
    \petit1 &   \jr  &   \jr  &   \jr  &   \jr  &   \jr  &   \je  \\
    \petit2 &   \+   &   \jr  &   \jr  &   \jr  &   \je  &        \\
    \petit3 &   \+   &   \+   &   \+   &   \je  &        &        \\
    \petit4 &   \jr  &   \jr  &   \je  &        &        &        \\
    \petit5 &   \jr  &   \je  &        &        &        &        \\
    \petit6 &   \je  &        &        &        &        &        \\
    \end{array}
  \]
  In this case, $P$ represents $\tilde\pi$ and this pipe dream~${\mathcal
  P}$ for $\tpi$ is reduced.
  \end{Example}
  Finally, we define the simplicial complex $\Gamma_\pi$ on $[n]^{2}$
  as follows. Note by construction that all crosses in each pipe dream
  for $L_\pi$ lie in the northwest $k\times n$ rectangle, and $k\leq
  n\leq m$. Then reverse the order of the first $k$ rows to obtain a
  pipe dream ${\hat \PP}$.  Then a face of $\Gamma_\pi$ is a subset of
  $[n]^2$ whose complement gives the positions of the crosses in
  ${\hat \PP}$.  In this way, $\Gamma_\pi$ and $L_\pi$ are isomorphic
  as simplicial complexes.  We will let $\mathit{Faces}(\Gamma_\pi)$
  denote the set of faces of $\Gamma_\pi$ (which we will identify with
  the corresponding pipe dreams).
  
  \begin{Example} \label{exa:upside_down}
  The crosses in the the first $k=3$ rows of the above pipe dream for
  $L_\pi$ given above are reversed to give the positions of the
  crosses of a pipe dream for $\Gamma_\pi$ as follows:
  %\begin{figure}[h]
  \[\begin{array}{ccccccc}
          &\perm1{}&\perm2{}&\perm3{}&\perm4{}&\perm5{}&\perm6{}\\
  \petit1 &     &     &     &     &     &     \\
  \petit2 & \+  &     &     &     &     &     \\
  \petit3 & \+  & \+  & \+  &     &     &     \\
  \petit4 &     &     &     &     &     &     \\
  \petit5 &     &     &     &     &     &     \\
  \petit6 &     &     &     &     &     &     \\
  \end{array}
  \ \ \mapsto \ \ \ 
  \begin{array}{ccccccc}
          &\perm1{}&\perm2{}&\perm3{}&\perm4{}&\perm5{}&\perm6{}\\
  \petit1 & \+  & \+  & \+  &     &     &     \\
  \petit2 & \+  &     &     &     &     &     \\
  \petit3 &     &     &     &     &     &     \\
  \petit4 &     &     &     &     &     &     \\
  \petit5 &     &     &     &     &     &     \\
  \petit6 &     &     &     &     &     &     \\
  \end{array}
  \]
  %\end{figure}
  The complement of positions of the crosses $\{(1,1), \ (1,2), \
  (1,3), \ (2,1)\}$ in the rightmost diagram give a face of
  $\Gamma_\pi$. Knowing the generators of $I_{31425}$,
  % (\ref{eqn:ess_gen3142}) and
  given Theorem~\ref{thm:grobnervex}, one checks that this complement
  is also a face in the Stanley--Reisner complex of $L_\pi$
  \comment{This had said $I_\pi$ but AK presumes that was a typo}.
  This agrees with the following theorem.
  \end{Example}
  
  \begin{Theorem}\label{thm:grobnervex'}
    Let $\pi\in S_n$ be a vexillary permutation.  With respect to any
    diagonal term order, the initial ideal of~$I_\pi$ is the
    Stanley--Reisner ideal for the simplicial complex~$\Gamma_\pi$.
    Moreover, $\Gamma_\pi$ and the subword complex~$L_\pi$ are
    isomorphic as simplicial complexes.
  \end{Theorem}
  \begin{proof}
  The last assertion follows from our discussions above.  For the
  remaining claim, if $\pi$ is Grassmannian, then by
  \cite[Theorem~B]{grobGeom}, $L_\pi$ is the Stanley--Reisner complex
  of the initial ideal of $I_\pi$ for any \emph{anti}diagonal term
  order.  Then it follows from Theorem~\ref{thm:grobnervex} and
  Lemma\verb=~\ref{lemma:grass_obs}= combined that the same is true of
  $\Gamma_\pi$ for $\init\,I_\pi$.
  
  Now suppose $\pi$ is not Grassmannian.  Then we can construct a
  sequence of vexillary permutations $\sigma_{1},\ \sigma_{2}, \ldots,
  \sigma_{t}=\pi$ satisfying the conclusions of
  Lemma~\ref{lemma:chain}.  We may assume by induction that
  $\Gamma_{\sigma_i}$ is the Stanley--Reisner ideal of
  $\init\,I_{\sigma_i}$ for $2\leq i<t$. Then by
  Proposition~\ref{p:inductionstep}, the partial degeneration
  $I_{\sigma_i}'$ of $I_{\sigma_i}$ satisfies
  \begin{equation}\label{eqn:partial_degen}
    I_{\sigma_i}'=I_{\sigma_{i+1}}\cap (I_{(\sigma_i)_P}+\<z_{l,m}\>)
    \hbox{\qquad\comment{$(l,m) \rightsquigarrow (p,q)$}}
  \end{equation}
  and this is a geometric vertex decomposition, where
  $(l,m)=(l_i,m_i)$ is an accessible box in $\Ess(\sigma_i)$.
  
  By fully degenerating both sides of
  % (\ref{eqn:partial_degen})
  the above equation to monomial ideals we obtain
  \[
    \init\,I_{\sigma_i}' = \init\,I_{\sigma_{i+1}} \cap
    \init(I_{(\sigma_i)_P}+\< z_{l,m}\>).
  \]
  Hence $\Gamma_{\sigma_i}$ is the union of the Stanley--Reisner
  complexes for the ideals $\init\,I_{\sigma_{i+1}}$ and
  $\init(I_{(\sigma_i)_P}+\< z_{l,m}\>)$.
  
  Since for any permutation $\rho$, the Schubert determinantal ideal
  $I_{\rho}$ is prime (see \cite{fulton:92} or \cite{grobGeom}), we
  know $I_{\sigma_{i+1}}$ is prime.  However, in addition,
  $I_{(\sigma_i)_P}+\< z_{l,m}\>$ is also prime: by
  Theorem~\ref{thm:grobnervex} the essential minors of
  $I_{(\sigma_i)_P}$ together with $z_{l,m}$ form a Gr\"obner basis of
  the ideal (since $z_{l,m}$ is relatively prime with these
  determinants).  Thus if a product of polynomials $fg$ is in
  $I_{(\sigma_i)_P}+\< z_{l,m}\>$ then $fg\mid_{z_{l,m}=0}$ is in
  $I_{(\sigma_i)_P}$. So we conclude that (say) $f\mid_{z_{l,m}=0}\in
  I_{(\sigma_i)_P}$.  Equivalently, $f\in I_{(\sigma_i)_P}+\<
  z_{l,m}\>$ (add back the terms involving $z_{l,m}$).
  
  By \cite[Theorem~1]{kalkbrener.sturmfels:initial}, the primality of
  the ideals $I_{\sigma_{i+1}}$ and $I_{(\sigma_i)_P}+\< z_{l,m}\>$
  shows that the two Stanley--Reisner complexes of their initial
  ideals are pure.  The squarefreeness of the generators shows that
  the initial ideals are also reduced, hence completely determined by
  their facets. Moreover, since the intersection of their monomial
  ideals is generated by the $\LCM$s of their respective generators,
  one can directly check that there are no facets in common
  \comment{check this!}.
  
  However, examining the facets of $\Gamma_{\sigma_i}$, those that do
  not have a cross $\textcross$ at position $(l,m)$ must be a facet of
  the Stanley--Reisner complex of $I_{\sigma_{i+1}}$. On the other
  hand, those that do have a cross there must be a facet of the
  Stanley--Reisner complex of $I_{(\sigma_i)_P}+\< z_{l,m}\>$. This
  allows us to conclude the induction step.
  \end{proof}
  %If $\pi$ is Grassmannian, then apply the reverse rows trick and
  %\cite{grobGeom} to conclude that $\Gamma_\pi$ is the desired
  %Stanley--Reisner ring as originally defined above. 
  %   
  %If $\pi$ is vexillary but not Grassmannian, then apply the argument
  %of the proof of the Gr\"obner basis theorem to reduce the problem to
  %the Grassmannian case. \comment{Details to come...}\end{proof}
%
}\end{excise}%

We could have waited to prove Theorem~\ref{thm:gvdshelling} until
after proving that the initial complex is a subword complex, and then
using the combinatorially proven fact from \cite{subword} that subword
complexes are shellable, instead of using geometric vertex decompositions. 
In \cite{Ksubword} we will reverse this argument, and use geometric vertex
decompositions to prove once more that subword complexes are shellable.

%%%%%%%%%%%%%%%%%%%%%%%%%%%%%%%%%%%%%%%%%%%%%%%%%%%%%%%%%%%%%%%%%%%%%%
\section{Flagged set-valued tableaux}\label{sec:fsvt}%%%%%%%%%%%%%%%%%
%%%%%%%%%%%%%%%%%%%%%%%%%%%%%%%%%%%%%%%%%%%%%%%%%%%%%%%%%%%%%%%%%%%%%%

\subsection{Set-valued tableaux and Grassmannian permutations}\label{sub:SVT}

\vspace{-.02ex}Formulae for certain Grothendieck polynomials
associated to a partition~$\lambda$ were given by Buch \cite{Buch}
(see Corollary~\ref{cor:buch_result}, below).  Naturally generalizing
the tableau formula for Schur polynomials, he gave his formula in
terms of \dfn{(semistandard) set-valued tableaux} with
shape~$\lambda$.  These are fillings $\tau: \lambda \to
PowerSet(\naturals)$ of the boxes in the shape of~$\lambda$ (in
English notation, the largest part along the top row) with nonempty
sets of natural numbers satisfying the following ``semistandardness''
conditions:
\begin{itemize}
\item
  if box $b \in \lambda$ lies above box $c \in \lambda$, then each
  element of $\tau(b)$ is strictly less than each element of
  $\tau(c)$; and
\item
  if box $b \in \lambda$ lies to the left of box $c \in \lambda$, then
  each element of $\tau(b)$ is less than or equal to each element of
  $\tau(c)$.
\end{itemize}
Let $\SVT(\lambda)$ denote the collection of all such tableaux
for a partition $\lambda$, as the nonsemistandard ones will be of
little interest in this paper. (They are more pertinent in \cite{KMYII}.)
For $\tau\in \SVT(\lambda)$, let $|\tau|$ denote the number of entries,
so $|\tau|\geq |\lambda|$ with equality exactly for ordinary (non-set-valued)
tableaux.

\begin{Example} \label{ex:tau}
The following is a set-valued tableau $\tau \in \SVT(\lambda)$ for the
partition $\lambda = (7,6,4,3,1)$:
\[
  \begin{array}{c}
  \tableau{1   & 1   & 1 & 1 & 1 & 1 & 1 \\
           2   & 2   & 2 & 2 & 2 & 2 \\
           3,4 & 4   & 4 & 4 \\ 
           5   & 5,6 & 6 \\
           6,7 & 7}
  \end{array}
=\
  \begin{array}{c}
  \tau.
  \end{array}
\]
\end{Example}

The goal of this section is to generalize and refine Buch's formula by
way of a Gr\"obner geometry explanation in terms of
Theorem~\ref{thm:grobnervex}.  The combinatorial \mbox{aspect} of this
story consists of a bijection between a certain set of pipe dreams and
a certain collection of set-valued tableau.  This generalizes the
bijection between (\mbox{ordinary}) semistandard Young tableaux and
reduced pipe dreams for Grassmannian permutations
(Proposition~\ref{p:SSYT}, below), which we learned from Kogan
\cite{KoganThesis} (on the other hand, see Example~\ref{exa:hard_bij}).

We begin with the bijection $\Omega$ from tableaux to pipe dreams.
Since its definition is just as easy to state for set-valued tableaux,
we work in that generality from the outset.
% Fix a Grassmannian permutation $\tpi \in S_N$ with descent at~$k$
% and associated partition~$\lambda = \lambda(\tpi)$.
Define $\Omega$ by associating to every $\tau \in \SVT(\lambda)$ a
% $k \times N$
pipe dream as follows:
% \[%\tag{$\Omega_{\,}$}
% \begin{array}[c]{c}
%   \hbox{\parbox[t]{.8\linewidth}{For each integer~$i$ that $\tau$
%   assigns to the box~$b$ of~$\lambda$, place a $\textcross$ in row~$i$
%   so that it lies on the diagonal containing the box~$b$.}}
% \end{array}
% \]
\begin{quote}
  For each integer~$i$ that $\tau$ assigns to the box~$b$,
%  of~$\lambda$, 
  place a $\textcross$ in row~$i$ so that it lies on the
  diagonal containing the box~$b$.
\end{quote}

\begin{Example} \label{ex:Omega}
The set-valued tableaux $\tau$ in Example~\ref{ex:tau} maps to
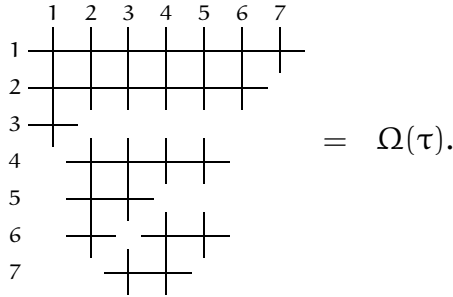
\[
\def\P#1{\perm{#1}{}}%
\def\peti#1{\petit{#1}}%
% \begin{array}{c}
%   \\\\\Omega(\tau)\ =\\\\
% \end{array}
%begin{array}{ccccccccccccccc}
%     &\P1&\P2&\P3&\P4&\P5&\P6&\P7&\P8&\P9&\P{10}&\P{11}&\P{12}&\P{13}&\P{14}\\
%peti1&\+ &\+ &\+ &\+ &\+ &\+ &\+ &   &   &      &      &      &      &      \\
%peti2&\+ &\+ &\+ &\+ &\+ &\+ &   &   &   &      &      &      &      &      \\
%peti3&\+ &   &   &   &   &   &   &   &   &      &      &      &      &      \\
%peti4&   &\+ &\+ &\+ &\+ &   &   &   &   &      &      &      &      &      \\
%peti5&   &\+ &\+ &   &   &   &   &   &   &      &      &      &      &      \\
%peti6&   &\+ &   &\+ &\+ &   &   &   &   &      &      &      &      &      \\
%peti7&   &   &\+ &\+ &   &   &   &   &   &      &      &      &      &      \\
%end{array}
\begin{array}{cccccccc}
      &\P1&\P2&\P3&\P4&\P5&\P6&\P7\\
\peti1&\+ &\+ &\+ &\+ &\+ &\+ &\+ \\
\peti2&\+ &\+ &\+ &\+ &\+ &\+ &   \\
\peti3&\+ &   &   &   &   &   &   \\
\peti4&   &\+ &\+ &\+ &\+ &   &   \\
\peti5&   &\+ &\+ &   &   &   &   \\
\peti6&   &\+ &   &\+ &\+ &   &   \\
\peti7&   &   &\+ &\+ &   &   &   
\end{array}
\ \ =\
  \begin{array}{c}
  \Omega(\tau).
  \end{array}
\]
\end{Example}
 
\begin{Proposition}\label{p:SSYT}
If~$\tpi \in S_N$ is a Grassmannian permutation with descent at~$k$
and associated partition~$\lambda = \lambda(\tpi)$, then $\Omega$
induces a bijection from the set of (ordinary) semistandard tableaux
of shape~$\lambda$ with entries at most~$k$ to the set of $k \times N$
reduced pipe dreams for~$\tpi$.  Moreover, in each $k \times N$
reduced pipe dream for~$\tpi$, the following hold.
\begin{alphlist}
\item
No pipe passes horizontally through one~$\textcross$ tile and
vertically through another.
% ~$\textcross$ tile.
\item
The row indices of the $\textcross$ tiles on the $i^\mathit{th}$
horizontal pipe from the top in $\Omega(\tau)$ are the values assigned
by~$\tau$ to the boxes in row~$i$ of~$\lambda$.
\item
The $\textcross$ tile in~$\Omega(\tau)$ corresponding to the box~$b$
of~$\lambda$ lies on the same diagonal as~$b$.
% itself.
\item
The box in row~$p$ and column~$q$ of~$\lambda$ corresponds to the
$\textcross$ in~$\Omega(\tau)$ at the intersection of the
$p^\mathrm{th}$ horizontal pipe from the top with the $q^\mathrm{th}$
vertical pipe from the~left.
% The $j^\mathit{th}$ vertical pipe from the left passes through the
% $j^\mathit{th}$ $\textcross$ tile on every horizontal pipe that has
% a $j^\mathit{th}$ $\textcross$ tile, and it passes through no other
% $\textcross$ tiles.
\end{alphlist}
% $\Omega$ maps $\FT$ bijectively to reduced pipe dreams; in the
% proof, do Grassmannian first and then the flagging for vexillary:
\end{Proposition}

\begin{Example} \label{ex:oltau}
Consider the ordinary semistandard Young tableaux $\ol\tau$ obtained
from the set-valued tableaux~$\tau$ in Example~\ref{ex:tau} by taking
only the smallest entry in each box.  Then, drawing the horizontal
pipes (with small bits of their crossings through vertical pipes to
make the picture easier to parse), $\Omega$ sends $\ol\tau$ to
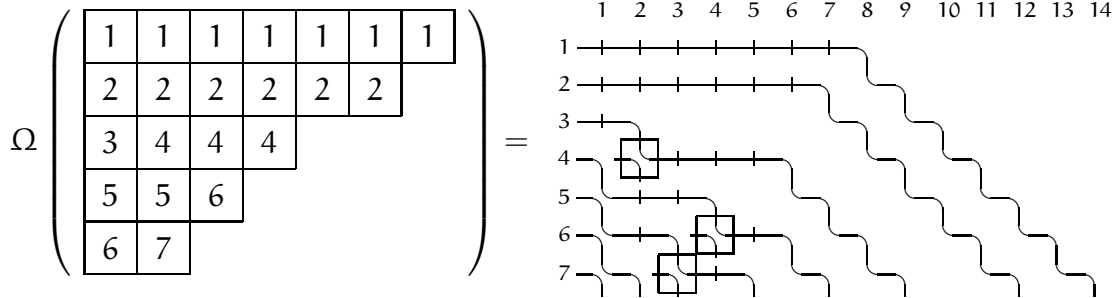
\[
\Omega\left(
\begin{array}{c}
\tableau{1   & 1   & 1 & 1 & 1 & 1 & 1 \\
         2   & 2   & 2 & 2 & 2 & 2 \\
         3   & 4   & 4 & 4 \\ 
         5   & 5   & 6 \\
         6   & 7}
\end{array}
\right)
=\
\def\P#1{\perm{#1}{}}%
\def\peti#1{\petit{#1}}%
\def\Be{\raisebox{-.5pt}{\makebox[0pt]{$\begin{array}{@{}|c|@{}}
\hline\phantom{\,\+\,}\\\hline\end{array}$}}\rj}%
\begin{array}{ccccccccccccccc}
      &\P1&\P2&\P3&\P4&\P5&\P6&\P7&\P8&\P9&\P{10}&\P{11}&\P{12}&\P{13}&\P{14}\\
\peti1&\ho&\ho&\ho&\ho&\ho&\ho&\ho&\ej&   &      &      &      &      &      \\
\peti2&\ho&\ho&\ho&\ho&\ho&\ho&\ej&\re&\ej&      &      &      &      &      \\
\peti3&\ho&\ej&   &   &   &   &\re&\ej&\re& \ej  &      &      &      &      \\
\peti4&\ej&\Be&\ho&\ho&\ho&\ej&   &\re&\ej& \re  & \ej  &      &      &      \\
\peti5&\rj&\ho&\ho&\ej&   &\re&\ej&   &\re& \ej  & \re  & \ej  &      &      \\
\peti6&\rj&\ho&\ej&\Be&\ho&\ej&\re&\ej&   & \re  & \ej  & \re  & \ej  &      \\
\peti7&\rj&\ej&\Be&\ho&\ej&\re&\ej&\re&\ej&      & \re  & \ej  & \re  & \ej  \\
\end{array}
\]
The three boxes containing elbow tiles will be explained in
Example~\ref{ex:abs}.
\end{Example}

\begin{proofof}{Proposition~\ref{p:SSYT}}
Let $\tau$ be a semistandard tableau of shape~$\lambda$ with entries
at most~$k$.  Construct a new tableau~$T$ (of shape $\lambda + \rho$,
where $\rho = (k,k-1,\ldots,2,1)$) by adding to the $i^\mathrm{th}$ row an
extra box filled with~$j$ for each $j = i, \ldots,k$ (and arranged to be
increasing along each row).  The map $\tau
\to \Omega(\tau)$ factors as $\tau \mapsto T \mapsto \Omega(\tau)$,
where each box $(p,q)$ of~$T$ filled with~$j$ corresponds to a tile at
$(j,q)$: the tile is a~$\textcross$ if $(p,q+1)$ is filled with~$j$
and a~$\textelbow$ otherwise.  The tiles in~$\Omega(\tau)$ not
assigned by~$T$ are all~$\textelbow$~tiles.

The semistandard condition on~$T$ guarantees that each $\textelbow$
arising from the last~$j$ in a row of~$T$ has another~$\textelbow$ due
south of it.  Hence there is a single pipe such that the numbers in
row~$i$ of~$T$ list the row indices of the tiles in~$\Omega(\tau)$
entered from the left by that pipe.  As every $\textcross$
in~$\Omega(\tau)$ arises from~$T$, this proves parts~(a) and~(b).
Part~(c)
% We apply \cite[Theorem~3.7]{BB}: Given a reduced pipe dream~$D$
% for~$\pi$, the Ferrers shape $\lambda(\pi)$ can be obtained from~$D$
% by repeatedly applying \emph{chute} operations (note that these are
% upside-down from those in~\cite{BB} because our elbow tiles
% are). Since every pipe passes only horizontally or only vertically
% through~$\textcross$ tiles, the only chutes available are those
% displacing a single $\textcross$ tile one unit diagonally
% northwest. The result follows by induction on the number of required
% chutes.%
follows easily from~(b) by considering each row of~$T$ separately and
using induction on the number of boxes in any fixed row.  Part~(d)
follows from part~(a), since the horizontal pipes can't cross one
another, and nor can the vertical~pipes.

It remains only to show that $\Omega$ induces the claimed bijection.
As the numbers of the indicated tableaux and reduced pipe dreams are
equal (both agree with the evaluation at $(1,\ldots,1)$ of the Schur
polynomial $s_\lambda(x_1,\ldots,x_k)$, which equals the Schubert
polynomial~$\SS_\tpi(x_1,\ldots,x_k)$), it is enough to show that
$\Omega(\tau)$ is a reduced pipe dream for~$\tpi$.  This follows
because $\Omega(\tau)$ has $|\lambda| = \mathrm{length}(\tpi)$
crossing tiles, and the pipe entering horizontally into row~$i$ exits
vertically out of row~$k$ through column $k-i+\lambda_i = \tpi(i)$, as
can be seen by counting its $\textcross$ tiles.%
\end{proofof}

Subword complexes are homeomorphic to balls or spheres, as shown in
\cite[Theorem~3.7]{subword}, where the interior and boundary faces
were characterized.  Here, the subword complex is $\Gamma_\tpi =
\Delta(Q,\tpi)$ for the full $k \times N$ rectangular word~$Q$.

\begin{Theorem} \label{t:SVT}
% \comment{EM: Here is the place to do the bijection \underline{in
% full} for Grassmannian permutations; the point is that the reader
% shouldn't be bogged down with the relatively trivial technicalities
% of flaggings before seeing really how to get from set-valued
% tableaux to nonreduced pipe dreams and back again}
If~$\tpi \in S_N$ is Grassmannian with descent at~$k$ and
partition~$\lambda = \lambda(\tpi)$, then $\Omega$ bijects the set
$\SVT_k(\lambda)$ of set-valued tableaux of shape~$\lambda$ with
entries at most~$k$ to the set of $k \times N$ pipe dreams~$P$
% such that $\tilde Q \minus P$ is an
whose elbow tiles form the vertex sets of interior
faces of~$\Gamma_\tpi$.
\end{Theorem}
\begin{proof}
The pipe dream $\Omega(\tau)$ for $\tau \in \SVT_k(\lambda)$ has
crosses in at most $k$ rows because the entries are at most~$k$, and
$\Omega(\tau)$ has at most~$N$ columns because $\lambda$ fits in a
rectangle of size $k \times (N-k)$.  No entry appears twice in~$\tau$
along any diagonal of~$\lambda$.  This together with semistandardness
implies that $\Omega$ is injective into its image.  What remains is to
show that the image of~$\Omega$ consists precisely of the interior
faces of~$\Gamma_\tpi$.

For ordinary tableaux~$\tau$ and facets~$\Omega(\tau)$
of~$\Gamma_\tpi$, this is Proposition~\ref{p:SSYT}.  For arbitrary
$\tau \in \SVT_k(\lambda)$, there is an associated ordinary
tableau~$\ol\tau$ obtained (as in Example~\ref{ex:oltau}) by taking
only the smallest entry in each box.  We will show that, similarly,
for each subword~$P$ of~$Q$ such that $Q \minus P$ is an interior face
of~$\Gamma_\tpi$, there is an associated reduced subword~$\ol P$
representing~$\tpi$.  To complete the proof, we will then construct
the set-valued tableau~$\tau$ satisfying $\Omega(\tau) = P$ starting
from the ordinary tableau~$\ol\tau$ satisfying $\Omega(\ol\tau) = \ol
P$.

For a reduced pipe dream~$\ol P$, say that a~$\textelbow$ tile in~$\ol
P$ is \dfn{absorbable} into~$\ol P$ if
% (i)~adding the~$\textcross$ to~$\ol P$ causes two of~$\ol P$'s pipes
% to cross a second time, and (ii)~the new~$\textcross$ is southeast
% of the original crossing of these two pipes in~$\ol P$.
the two pipes passing through it intersect in a~$\textcross$ tile to
its northwest (see Example~\ref{ex:abs}, below).  It is immediate from
the definition that a tile is absorbable if and only if the
corresponding reflection in $Q \minus \ol P$ is absorbable in the
sense of \cite[Section~4]{subword}.  Therefore it follows from
\cite[Theorem~3.7]{subword} that a pipe dream~$P$ is the complement of
an interior face $Q \minus P$ of~$\Gamma_\tpi$ if and only if $P$ is
obtained from a reduced pipe dream~$\ol P$ for~$\tpi$ by changing some
of its absorbable $\textelbow$ tiles into $\textcross$ at will.

Proposition~\ref{p:SSYT}(a) allows us to distinguish between
\emph{horizontal} and \emph{vertical} pipes in any reduced pipe
dream~$\ol P$ for~$\tpi$.  We claim that if a horizontal and vertical
pipe cross at~$\textcross$ and pass through a $\textelbow$ tile
southeast of it, then the $\textelbow$ lies on the same diagonal as
the~$\textcross$ and occurs in a row that is strictly north of the
next $\textcross$ down (if there is one) on the vertical pipe.  This
suffices because altering such~$\textelbow$ tiles to~$\textcross$
tiles corresponds, by defintion of~$\Omega$, to inserting extra
entries in the box of~$\ol\tau$ corresponding to the
original~$\textcross$ tile of~$\ol P$.  The claim holds because once
the vertical pipe passes downward through the next horizontal pipe,
that next horizontal pipe separates the vertical pipe
% it can never again approach the
from the original horizontal pipe.%
\end{proof}

\begin{Example} \label{ex:abs}
The shape $\lambda$ in Examples~\ref{ex:tau} and~\ref{ex:oltau} is
$\lambda(\tpi)$ for a Grassmannian permutation $\tpi \in S_{N=14}$
% for $N = 14$
with descent at $k = 7$.  The absorbable tiles in
Example~\ref{ex:oltau} lie in the boxes at $(4,2)$, $(6,4)$, and
$(7,3)$.  Altering these to $\textcross$ tiles yields the nonreduced
pipe dream in Example~\ref{ex:Omega}; as in the proof of
Theorem~\ref{t:SVT}, these absorbable tiles also correspond to the
extra entries needed to make the set-valued tableau in
Example~\ref{ex:tau} from the ordinary tableau in
Example~\ref{ex:oltau}.
\end{Example}

%       For example, originally, 
%Buch \cite{Buch} generalized the tableaux formula for
%Schur polynomials by giving formula for the generating series
%$\GG_{\lambda}(\xx)$ for all tableaux in $\SVT(\lambda)$
%(say) whose
%entries are from the set $\{1,2,\ldots,M\}$.
%For example, if $\lambda=(2,1)$ and $M=2$, the reader can check that
% there are three set-valued tableaux and hence
%\[\GG_{\lambda}(x_1,x_2)=\sum_{t\in \SVT(2,1)}\prod_{b\in
%t}{\wt_\xx(b)}
%=x_1^2 x_2 + x_1 x_2^2 - x_1^2 x_2^2.\]
%Other formulae for $\GG_{\lambda}(\xx)$ were previously known  
%\cite{fomin.stanley,FKyangBax}: for
%any permutation $\pi$ there is a 
%{\bf Grothendieck polynomial} $\GG_\pi(\xx,\yy)$;
%$\GG_{\lambda}(\xx)$ is $\GG_\pi(\xx,{\bf 0})$ for
%a permutation $\pi$ with $\lambda=\lambda(\pi)$. 
%We will mainly need the interpretation of $\GG_\pi(\xx,\yy)$
%and their lowest degree homogeneous component, the {\bf Schubert
%polynomial} $\SS_\pi(\xx,\yy)$) as a
%$K$-polynomial and multidegree of the matrix Schubert variety
%$\ol X_\pi$, 
%as proved in \cite{grobGeom} (we also refer the reader to this paper 
%and the references therein for earlier aspects of these polynomials).
% 
%       We will need to introduce a class of set-valued tableaux
%with row bounds on the entries, depending on a vexillary permutation
%$\pi\in S_n$.
% 
%       There is a natural correspondence from the set of outside corners of 
%$\lambda(\pi)$ to the accessible boxes of $D(\pi)$, which 
%takes the rightmost box $c$ from the lowest row in $\lambda$ with $b$ boxes
%to the rightmost box in the lowest row of $D(\pi)$ with $b$ boxes.

\subsection{Flaggings}\label{sub:flaggings}

Given a partition $\lambda$, a \dfn{flagging} $f$ of~$\lambda$ is a
natural number assigned to each row of $\lambda$.  Suppose that
$\lambda = \lambda(\pi)$ for some vexillary permutation~$\pi$, and
recall the partition~$\mu(\pi)$ from Section~\ref{section:subword}.
The \dfn{flag} $f_\pi$
% is defined as follows:
% if the $i^\mathrm{th}$ row of $D(\pi)$ is nonempty, let $e_i$ be the
% greatest integer $j\geq i$ such that the $j^\mathrm{th}$ row has at
% least as many boxes of $D(\pi)$ as the $i^\mathrm{th}$; then $f_\pi$
% is the sequence of integers $e_i$ sorted into increasing order.
% This flag $f_\pi$ provides a flagging of $\lambda$, by assigning to
% the row $i$ of the partition.  It is known that from $\lambda(\pi)$
% and $f_\pi$ one can uniquely reconstruct $\pi$, see,
% e.g.,~\cite[Proposition~2.2.10]{Manivel}.
assigns to row~$i$ of~$\lambda$ the row index of the southeastern box
of $\mu(\pi)$ that lies on the same diagonal as the last box
$(i,\lambda_i)$ in row~$i$ of~$\lambda$.

If $f$ is a flagging of~$\lambda$, call a set-valued tableau~$\tau$
\dfn{flagged by~$f$} if each box~$b$ in the $i^\mathrm{th}$ row
of~$\lambda$ satisfies $\max\tau(b) \leq f(i)$.  In particular,
if~$\tau$ is an ordinary tableau (each set a singleton), this
definition of flagged tableaux is the usual one, in \cite{Wachs}, for
example.  Let $\FST(\pi)$ denote the collection of set-valued tableaux
of shape~$\lambda(\pi)$ flagged by~$f_\pi$.  Also, let $\FT(\pi)$
denote the subset consisting of flagged semistandard (ordinary)
tableaux, in which, by definition, all the sets are singletons.

\begin{Example}
If $\pi$ is the vexillary permutation from
Example~\ref{example:pi_CP}, then $\lambda(\pi)$ is the partition from
Example~\ref{ex:tau}, and the tableau~$\tau$ in Example~\ref{ex:tau}
obeys the flagging $f_\pi=(1,2,4,6,7)$; that is, $\tau$ lies in
$\FST(\pi)$.  In general, the flag of a vexillary permutation need not
list the indices of the nonempty rows of~$D(\pi)$; indeed, for the
permutation~$\pi_C$ in Example~\ref{example:pi_CP}, the flagging is
$f_{\pi_C} = (1,2,4,6,6)$.%
\end{Example}

Let $\xx=(x_1,x_2,\ldots)$ and $\yy=(y_1,y_2,\ldots)$ be two
collections of commuting indeterminates.  
To any permutation $\pi\in S_n$, Lascoux and Sch\"utzenberger
associated a \dfn{(double) Grothendieck polynomial} $\GG_\pi(\xx,\yy)$
\cite{LS82}.  Our convention here is that $\GG_\pi(\xx,\yy)$ means the
same thing as in \cite{grobGeom}, which would be called
$\GG_\pi(\xx,\yy^{-1})$ in \cite{subword}, and is obtained from the
polynomial called $\mathfrak
L\begin{array}{@{}l@{}}\\[-4ex]\scriptstyle (-1)\\[-1.5ex]
\scriptstyle\, \pi\\[-1ex]\end{array}(y,x)$ in \cite{FK94} by
replacing each $x$ and~$y$ variable there with $1-x$ and
$\frac{1}{1-y}$, respectively.  Let $\SS_\pi(\xx,\yy)$ denote the
\dfn{(double) Schubert polynomial}, which is the lowest homogeneous
degree component of $\GG_\pi(\mathbf{1}-\xx,\mathbf{1}-\yy)$ when this
rational function is expressed as a series in positive powers of $\xx$
and~$\yy$.  These double Schubert polynomials are the same as those in
\cite{grobGeom}.  Combinatorial formulae for Schubert and Grothendieck
polynomials are known in terms of pipe dreams \cite{FKyangBax}, and
Gr\"obner geometry explanations of these formulae were given
in~\cite{grobGeom,subword}.  However, the approach given there does
not explain the tab\-leau formulae given here, which are different and
only apply when $\pi$ is vexillary.  Hence we now give a Gr\"obner
geometry explanation of these tableau formulae.

\begin{excise}{
  We will be interested in generating functions for tableaux in $\xx$
  and~$\yy$ with respect to the following weights: for $\tau \in
  \SVT(\lambda)$ and a box $b \in \lambda$,~define
% \[
%   \wt_\xx(b)=\prod_{i\in t(b)}(-1)^{\mid t(b)\mid-1} x_i
%   \qquad\hbox{and}\qquad\wt_\yy(b)=\prod_{i\in t(b)}y_{i-r(b)+c(b)},
% \]
\[
%  \wt^\tau_\xx(b)= -\prod_{i\in \tau(b)} (-x_i) \qquad\hbox{and}\qquad
%  \wt^\tau_{\xx,\yy}(b)=-\prod_{i\in\tau(b)}(x_iy_{i+j(b)}-x_i-y_{i+j(b)}),
  \wt^\tau_{\xx,\yy}(b) =
  -\prod_{i\in\tau(b)}\Big(\frac{x_i}{y_{i+j(b)}} - 1\Big),
  \]
  where $j(b) = r(b) - c(b)$ is the difference of the row and column
  indices of $b \in \lambda$.
}\end{excise}

\begin{Theorem} \label{t:FSVT}
If $\pi\in S_n$ is a vexillary permutation, then $\Omega$ bijects the
facets and interior faces of~$\Gamma_\pi$ with $\FT(\pi)$
and~$\FST(\pi)$, respectively.  Consequently, we have the following
formulae for the double Schubert and double Grothendieck polynomials
associated to~$\pi$:
\begin{eqnarray*}
  \SS_\pi(\xx,\yy) &=& \ \sum_{\tau\in \FT(\pi)\,} \prod_{e\in \tau} \ 
  (x_{\val(e)} - y_{\val(e) + j(e)})
\\
  \GG_\pi(\xx,\yy) &=& \sum_{\tau\in \FST(\pi)\,} (-1)^{|\tau|-|\lambda|}
  \prod_{e\in \tau} 
  \Big(1 - \frac{x_{\val(e)}}{y_{\val(e)+j(e)}} \Big)
\end{eqnarray*}
where each product is over each entry $e$ of $\tau$, whose numerical
value is denoted $\val(e)$, and
where $j(e) = c(e) - r(e)$ is the difference of the row and column
indices. The sign $(-1)^{|\tau|-|\lambda|}$ alternates with the number of
``excess'' entries in the set-valued tableau.

%\begin{eqnarray*}
%  \SS_\pi(\xx,\yy) &=& \ \sum_{\tau\in \FT(\pi)\,} \prod_{\ (p,q)\in
%  \lambda(\pi)} (x_i - y_{i + p - q})
%\\
%  \GG_\pi(\xx,\yy) &=& \sum_{\tau\in \FST(\pi)\,} \prod_{\,b\in
%  \lambda(\pi)} \wt^\tau_{\xx,\yy}(b)
%\end{eqnarray*}
\end{Theorem}
\begin{proof}
\begin{excise}{%
  \comment{EM: What follows can't be a proof, because the ``interior''
  condition doesn't play a role:} 
  
          It remains to prove that $\Omega(\FST(\pi))$ is equal to
  the set of pipe dreams that represent faces of $\Gamma_\pi$.  For
  $i \geq 0$ let $\Omega(\FST)_i$ denote the subset of
  $\Omega(\FST(\pi))$ with $i+|\lambda(\pi)|$ crosses.  We prove the
  desired claim by showing by induction on $i$ that
  $\Omega(\FST(\pi)_i)$ is equal to the subset
  $\mathit{Faces}(\Gamma_\pi)_i$ of
  %  set $\mathit{Faces}(\Gamma_\pi)$
  the set of faces of~$\Gamma_\pi$ with the same number of crosses.
  The base case of our induction is Proposition~\ref{prop:basecase}.
  
  %       We first prove that the map provides a bijection between 
  %$\FT(\pi)\subset \FST(\pi)$ the facets of $\Gamma_\pi$.  It follows
  %from \comment{reference? due to Kogan?} that the pipe dreams for 
  %the Grassmannian permutation $\tpi$ which correspond to facets of 
  %$L_\pi$ are precisely those that under the ``upside down'' bijection 
  %of Theorem~\ref{thm:grobnervex} 
  %(illustrated in Example~\ref{exa:upside_down}) we obtain precisely
  %the pipe dreams in the image of $\Omega$.  Then by
  %Theorem~\ref{thm:grobnervex}
  %these pipe dreams are the facets of $\Gamma_\pi$.  Hence we conclude
  %that $\Omega$ is a bijection between $\FT(\pi)$ and the facets of
  %$\Gamma_\pi$
  %follows.
  
          For $i > 0$ we next show that 
  $\Omega(\FT(\pi))_i \subseteq \mathit{Faces}(\Gamma_\pi)_i$.
  Let ${\hat \PP} = \Omega(\tau) \in \Omega(\FT(\pi))_i$ and
  pick any box of $t$ which has at least two entries.  Then by
  induction, after removing the cross in ${\hat \PP}$ that corresponds
  to the 
  largest entry of that box we obtain a pipe dream ${\hat \PP}'$
  in $\mathit{Faces}(\Gamma_\pi)_{i-1}=\Omega(\FT(\pi))_{i-1}$.
  Turning this latter pipe dream upside
  down, we get $\PP'\in L_\pi$ which therefore
  corresponds to a subword $P'$ of $Q$
  that contains $\tpi$.  It is then not hard to check that 
  the semistandard conditions on $t$ imply that
  after replacing the cross we removed from ${\hat \PP}$ 
  (and correspondingly adding a cross to $\PP'$) 
  we see that the subword $P$ (corresponding to $\PP$) 
  still contains $\tpi$.
  Hence the desired containment holds.
  
          It remains to check the other containment.  Let 
  ${\hat \PP}\in \mathit{Faces}(\Gamma_\pi)_i$ for some $i>0$ and
  let $\PP$ be the corresponding pipe dream in $L_\pi$.
  Define a cross in $\PP$ 
  to be \dfn{removable} if after removing that cross, we remain in
  $L_\pi$. We will need:
  \comment{EM: This definition of removable is not what is required;
  one needs that removing the cross in question doesn't change the
  Demazure product, which is touchy to define directly; instead, I've
  resorted to the notion of crosses being ``absorbable'' in the sense
  of \cite[Section~4]{subword}, which is easier to define pictorially}
  \begin{Lemma} 
  \label{lemma:removable}
  For any vexillary permutation $\pi\in S_n$ the following
  hold:
  \begin{itemize}
  \item[(a)] any non facet of $L_\pi$ contains a removable cross; and
  \item[(b)] a cross in a pipe dream for $L_\pi$ is removable if and
  only if the two strands that intersect at that cross also intersect
  at another cross to the southwest.  Moreover the two crosses appear
  in the same antidiagonal. 
  \end{itemize}
  \end{Lemma}
  \begin{proof}
  Parts (a) and the first part of (b) follow easily from the results
  of \cite{subword} (see in particular,
  the comments after the proof of Corollary~5.5).  For part (b), 
  since the pipe dreams in question are for a Grassmannian
  permutation, the latter
  assertion follows from the results of \cite{BB}, or e.g., 
  Proposition~\ref{prop:basecase}.
  \end{proof}
  
          By Lemma~\ref{lemma:removable}(a) we can find
  and delete a removable cross of $\PP$ 
  to obtain a one smaller pipe dream
  $\PP'\in L_\pi$.  By induction, $\PP'=\Omega(\tau')$ 
  for some set-valued tableau $\tau'$.  Then
  Lemma~\ref{lemma:removable}(b) together with the semistandard
  conditions on $\tau'$ implies that replacing the cross we removed
  corresponds to adding an entry to a box of $\tau'$ while
  remaining in $\FST(\pi)$.  Hence we conclude 
  $\Omega(\FT(\pi))_i \supseteq \mathit{Faces}(\Gamma_\pi)_i$.
  
          Lastly, the above bijection together with
  \cite[\comment{which theorem?}]{subword} proves the displayed
  formulae for $\SS_\pi(\xx, \yy)$ and $\GG_\pi(\xx,\yy)$
}\end{excise}%
Let $\tpi$ be the Grassmannian permutation of descent~$k$ associated
to~$\pi$ by Lemma~\ref{lemma:chain}.  The image under $\Omega$
of~$\FST(\pi)$ consists exactly of the pipe dreams in $\Omega(\SVT_k)$
whose $\textcross$ tiles all lie inside the Ferrers shape~$\mu(\pi)$.
Therefore the first sentence follows immediately from
\cite[Theorem~3.7]{subword} and Theorem~\ref{t:SVT}.

Consider the $\ZZ^{2n}$-grading on~$\kk[\zz]$ from
Section~\ref{sub:schub}.  The desired double Grothen\-dieck polynomial
equals the $\ZZ^{2n}$-graded $\mathit{K}$-polynomial of the quotient
$\kk[\zz]/I_\pi$ by \cite[Theorem~A]{grobGeom}.  This statement holds
for $\kk[\zz]/\IN\,I_\pi$, too, since $\mathit{K}$-poly\-nomials are
invariant under taking initial ideals \cite[Theorem~8.36]{cca}.  On
the other hand, the $\ZZ^{2n}$-graded $\mathit{K}$-polynomial
of~$\kk[\zz]/\IN\,I_\pi$ is obtained from its $\ZZ^{n^2}$-graded
% $\mathit{K}$-polynomial
counterpart by replacing~$z_{pq}$ with~$x_p/y_q$.  Using the
definition of~$\Gamma_\pi$ as a subword complex, one calculates the
$\ZZ^{n^2}$-graded $\mathit{K}$-polynomial of~$\kk[\zz]/\IN\,I_\pi$ as
in \cite[Theorem~4.1]{subword}: it is, by \cite[Theorem~3.7]{subword},
the sum over the pipe dreams~$P$ corresponding to interior faces $Q
\minus P$ of~$\Gamma_\pi$ of products $(-1)^{|P|-\ell}\prod_{(p,q) \in
P} (1-z_{pq})$, where $\ell = |D(\pi)|$ is the length of~$\pi$.  The
formula for~$\GG_\pi(\xx,\yy)$ results because, by definition
of~$\Omega$, the summand $\prod_{b \in
\lambda(\pi)}\wt^\tau_{\xx,\yy}(b)$ for $\tau \in \FST(\pi)$ equals
the appropriately signed product of factors $(1-x_p/y_q)$ over all
$(p,q)$ such that $\Omega(\tau)$ has a~$\textcross$ tile at~$(p,q)$.

Substituting $1-x_p$ for $x_p$ and $1-y_q$ for $y_q$ in
$\wt^\tau_{\xx,\yy}(b)$ results in a power series whose lowest term
has degree equal to the cardinality of~$\tau(b)$.  Taking the lowest
degree terms in~$\GG_\pi(\mathbf{1}-\xx,\mathbf{1}-\yy)$ therefore
yields a sum over honest (that is, not set-valued) tableaux, and the
formula for~$\SS_\pi(\xx,\yy)$ follows.%
\end{proof}

In particular, we obtain a Gr\"obner geometry %/Gr\"obner degeneration
explanation of the following result due to Buch
\cite[Theorem~3.1]{Buch}.

\begin{Corollary}[\cite{Buch}] \label{cor:buch_result}
If $\pi$ is a Grassmannian permutation and $\lambda=\lambda(\pi)$,
then
\[
  \GG_{\lambda}(\mathbf{1}-\xx) = \sum_{\tau\in \SVT(\lambda)\,}
% \prod_{b \in \tau}\wt^\tau_\xx(b),
  (-1)^{|\tau|-|\lambda|} \prod_{e\in \tau(b)} x_{\val(e)},
\]
where $\GG_\lambda(\mathbf{1}-\xx)$ is obtained from
$\GG_\lambda(\xx,\yy)$ by replacing each $x_p$ and~$y_q$ with $1-x_p$
and~$1$.
\end{Corollary}

\begin{Example} \label{exa:hard_bij}
We emphasize that our proof of Theorem~\ref{t:FSVT} is a consequence
of the Gr\"obner geometry, but it is \emph{not}\/ a direct
combinatorial consequence of the Fomin--Kirillov formulae \cite{FK94,
FKyangBax}.  We leave it as a challenge to find such an explanation,
even in the case of double Schubert polynomials.  For example, let
$\pi=\left(\begin{array}{cccc} 1 & 2 & 3 & 4 \\ 1 & 4 & 3 & 2 \\
\end{array}\right)$.  Then $\lambda(\pi)=(2,1)$ and $f_{\pi}=(2,3)$.
Theorem~\ref{t:FSVT} computes
\begin{eqnarray*}
\SS_\pi(\xx,\yy) 
&=& (x_2-y_2)(x_2-y_3)(x_3-y_2)+ (x_1-y_1)(x_2-y_3)(x_3-y_2)\\
& & \mbox{}+(x_1-y_1)(x_1-y_2)(x_3-y_2) + (x_1-y_1)(x_2-y_1)(x_2-y_3)\\
& & \mbox{}+(x_1-y_1)(x_1-y_2)(x_2-y_1).
\end{eqnarray*}
On the other hand, using the formula of \cite{FKyangBax} gives
\begin{eqnarray*}
\SS_\pi(\xx,\yy) 
&=& (x_1-y_3)(x_2-y_1)(x_3-y_1) + (x_1-y_2)(x_1-y_3)(x_3-y_1)\\
& & \mbox{}+(x_2-y_1)(x_2-y_2)(x_3-y_1) + (x_1-y_2)(x_2-y_1)(x_2-y_2)\\
& & \mbox{}+(x_1-y_2)(x_1-y_3)(x_2-y_2).
\end{eqnarray*} 
One approach to relating these two formulae would be to prove that
after any permutation of the rows, the essential minors remain a
Gr\"obner basis for any diagonal term order.  Then each permutation of
the rows would give a different formula, and one might be able to
relate the formulae associated to permutations that are adjacent in
Bruhat order.  The two formulae above correspond to the identity and
long-word permutations of the rows.
\end{Example}

%%%%%%%%%%%%%%%%%%%%%%%%%%%%%%%%%%%%%%%%%%%%%%%%%%%%%%%%%%%%%%%%%%%%%%
\section{The diagonal Gr\"obner basis theorem for Schubert ideals is sharp}
%%%%%%%%%%%%%%%%%%%%%%%%%%%%%%%%%%%%%%%%%%%%%%%%%%%%%%%%%%%%%%%%%%%%%%
\label{section:sharp}

Our goal in this section is to prove the converse of
Theorem~\ref{thm:grobnervex}.

Let $A_\pi$ denote the union over all $(p,q) \in \Ess(\pi)$ of the
sets of minors of size $1+\rr pq\pi$ in the northwest $p\times q$
corner $\zz_{p \times q}$ of the generic matrix~$\zz$.  Define $B_\pi$
similarly, except take the union over all $(p,q)$ in the $n \times n$
grid.  Both sets generate~$I_\pi$, as shown in \cite{fulton:92} (see
also \cite[Chapter~15]{cca} for an exposition).
% The following result explains why our Gr\"obner basis theorem does
% not hold for the Schubert determinantal ideal of any non-vexillary
% permutation.

\begin{Theorem} \label{thm:neg}
If $\pi\in S_n$ is a permutation that is not vexillary, then neither
$A_\pi$ nor $B_\pi$ is a Gr\"obner basis of $I_\pi$ under any diagonal
term order.
\end{Theorem}

In what follows,
% we will use the following terminology and notation.
we say that a pipe dream~$P$ \dfn{poisons} a diagonal or a minor if at
least one member of the diagonal (or the diagonal term of the minor)
coincides with a~$\textcross$ in~$P$.  We will also say that $P$ is a
\dfn{poisoning of a set} of diagonals or minors if it poisons each
element in the set.  The poisoning of a set is \dfn{minimal} if by
removing any cross from~$P$, some element in the set is no longer
poisoned.  Recall the diagram~$D(\pi)$ from
Section~\ref{sub:matrixSchub}, and let $C(\pi)$ denote the $n \times
n$ pipe dream formed by placing $\textcross$ tiles in the boxes
of~$D(\pi)$.

Our proof of Theorem~\ref{thm:neg} is based on the following.

\begin{Proposition} \label{prop:poison}
The cross diagram $C(\pi)$ is a poisoning of~$A_\pi$.  Moreover, $C(\pi)$
is a minimal poisoning of $A_\pi$ if and only if $\pi$ is vexillary.
\end{Proposition}

The proof of this proposition requires a number of intermediate
results.  As in Section~\ref{sec:gvdvmSv}, identify each
permutation~$\pi$ with its dot-matrix.

\begin{Lemma} \label{lemma:main_obs_unpois}
For any diagonal~$\delta$ not poisoned by $C(\pi)$, the following
holds.
\begin{alphlist}
\item
Each element of~$\delta$ is either south of a dot of~$\pi$ in the same
column, or east of a dot of~$\pi$ in the same row.
\item
No dot of $\pi$ has an element of~$\delta$ south of it in the same
column and a different element of~$\delta$ east of it in the same row.
\end{alphlist}
\end{Lemma}
%The diagram $D(\pi)$ naturally splits into \dfn{connected components},
%each of which has a unique northwest corner.
\begin{proof}
These statements are immediate consequences of the definitions.
\end{proof}

We identify diagonals with the products of the corresponding variables
in the array~$\zz$.  Recall that $\zz_{p \times q}$ denotes the
northwest $p \times q$ submatrix of~$\zz$.

\begin{Corollary} \label{cor:unpois_diag}
The following hold for any $(p,q) \in D(\pi)$.
\begin{alphlist}
\item
Any diagonal in $\zz\sub pq$ that is not poisoned by~$C(\pi)$ has size
at most~$\rr pq\pi$.
\item
If $(p,q)$ is any maximally northwest box in~$D(\pi)$ with the
property that at least two of the $\rr pq\pi = \rr{p-1,}{q-1}\pi$ dots
in $\pi\sub{(p-1)}{(q-1)}$ do not lie on a diagonal, then any diagonal
in $\zz\sub{(p-1)}{(q-1)}$ not poisoned by~$C(\pi)$ has size at most
$\rr pq\pi - 1$.
\end{alphlist}
\end{Corollary}
\begin{proof}
Part (a) follows immediately from the definition of $\rrr\pi$ and
Lemma~\ref{lemma:main_obs_unpois}.

Our hypotheses in part~(b) imply that some pair of dots
in~$\pi\sub{(p-1)}{(q-1)}$, say at $(i,j)$ and~$(k,\ell)$, forms an
antidiagonal (where, say $(i,j)$ occurs to the southwest of $(k,\ell)$).  
Additionally, in~$\zz\sub{(p-1)}{(q-1)}$, every column
from $\ell$ through~\mbox{$q-1$} has a dot of~$\pi$: if not, some box
of~$D(\pi)$ in row~$p$ strictly west of $(p,q)$ would contradict our
assumptions on~\mbox{$(p,q)$}.  Similarly, every row from $i$
through~\mbox{$p-1$} has a dot of~$\pi$.

% Now suppose that $\pi\sub{(p-1)}{(\ell-1)}$ has $s$ dots.  By our
% choice of~$\ell$, all of these dots form a diagonal, and it follows
% from Lemma~\ref{lemma:main_obs_unpois} that any diagonal in
% $\zz\sub{(p-1)}{(\ell-1)}$ not poisoned by~$C(\pi)$ has size at
% most~$s$.  By construction, the number of columns from~$\ell$ to
% $q-1$ is $\rr pq\pi - s$.  Therefore, since $\delta$ has maximal
% size, it must restrict to a diagonal of size~$s$ in
% $\zz\sub{(p-1)}{(\ell-1)}$ and it must intersect each of the columns
% from~$\ell$ through~\mbox{$q-1$}.  These two observations show that
% the 
Assume that some diagonal~$\delta$ in $\zz\sub{(p-1)}{(q-1)}$ not
poisoned by~$C(\pi)$ has maximal length~$\rr pq\pi$.  Suppose that for
some~$h$ satisfying $\ell \leq h \leq q-1$, no element of~$\delta$
lies in column~$h$ weakly south of the dot of~$\pi$ there, and
choose~$h$ maximal with this property.  It follows that
if~$\,\bullet\,$ is the dot of~$\pi$ in column~$h$, then any element
of~$\delta$ in the hook of~$\,\bullet\,$ lies in the hook of some
other dot of~$\pi$ (by maximality of $h$), 
and $\delta$ is forced to have size less than $\rr
pq\pi$ by Lemma~\ref{lemma:main_obs_unpois}.  Therefore, every column
from~$\ell$ through $q-1$ has an element of~$\delta$ weakly south of
the dot of~$\pi$ there.  Similarly, every row from~$i$ through $p-1$
has an element of~$\delta$ weakly east of the dot of~$\pi$~there.
% The intersection of~$\delta$ with the rectangle
% in~$\zz\sub{(p-1)}{(q-1)}$ whose northwest corner is~$(i,\ell)$ is a
% diagonal of size at most $\min(p-i,q-\ell)$, and this minimum is at
% least~$1$.

Where is the element of~$\delta$ in column~$\ell$?  If it is weakly
north of row~$i$, then any element of~$\delta$ east of~$(i,j)$ in
row~$i$ also lies south of some other dot of~$\pi$ (in a column
from~$\ell$ to $q-1$).  On the other hand, if the element of~$\delta$
in column~$\ell$ is weakly south of row~$i$, then it lies east of some
other dot of~$\pi$ (in a row from~$i$ to $p-1$).  In both cases,
$\delta$ is forced to have size less than $\rr pq\pi$ by
Lemma~\ref{lemma:main_obs_unpois}, a contradiction.%
\end{proof}

\begin{proofof}{Proposition~\ref{prop:poison}}
Part~(a) of Corollary~\ref{cor:unpois_diag} proves the first
assertion.  For the remaining assertion, suppose first that $\pi$ is
vexillary.  Pick any $(p,q) \in \Ess(\pi)$.
Corollary~\ref{cor:vex_equiv} shows that the diagonal formed by the
dots of~$\pi\sub{(p-1)}{(q-1)}$ together with $(p,q)$ form a diagonal
of size~$1 + \rr pq\pi$.  Hence we cannot remove $(p,q)$ from~$C(\pi)$
and remain a poisoning.  Thus $C(\pi)$ is a minimal poisoning.

Conversely, suppose that $\pi$ is not vexillary.  By
Corollary~\ref{cor:vex_equiv}, we find can find $(p,q)\in D(\pi)$ that
is maximally northwest with the property that at least two dots
of~$\pi\sub{(r-1)}{(s-1)}$ are in antidiagonal position.  We claim
that removing the cross of~$C(\pi)$ from $(p,q)$ results in a smaller
poisoning.  If this is not true, we can find $(i,j)\in \Ess(\pi)$ to
the southeast of $(p,q)$ and a diagonal of size $1 + \rr ij\pi$ in
$\zz\sub ij$ unpoisoned by $C(\pi) \minus (p,q)$ but containing
$(p,q)$.  By Corollary~\ref{cor:unpois_diag}(b) the part of this
diagonal in~$\zz\sub pq$ (including $(p,q)$ itself) has size at
most~$\rr pq\pi$.  The rest of this diagonal, which lies strictly
south and strictly east of $(p,q)$, has size at most $\rr ij\pi - \rr
pq\pi$ by Lemma~\ref{lemma:main_obs_unpois}, because this rank
difference equals the number of dots of~$\pi\sub ij$ in the union of
all rows~\mbox{$> p$} and columns~\mbox{$> q$}.  Hence the maximal
size of our unpoisoned diagonal is~$\rr ij\pi$, a contradiction.  This
completes the proof of Proposition~\ref{prop:poison}.%
\end{proofof}

For the proof of Theorem~\ref{thm:neg} we will also need the
following.

\begin{Proposition} \label{prop:divisible}
For any permutation~$\pi$, the diagonal term of each minor in~$B_\pi$
is divisible by the diagonal term of some minor in~$A_\pi$.
\end{Proposition}
\begin{proof}
Fix a minor in $B_\pi$ of size $1 + \rr ij\pi$ in $\zz\sub ij$.  If
$(i,j) \in D(\pi)$ then the result follows by definition, so assume
otherwise.

As the minor's diagonal~$\delta$ is larger than the number of dots in
$\pi\sub ij$, Lemma~\ref{lemma:main_obs_unpois} implies that some
element in~$\delta$ must lie in~$D(\pi)$.  Let $(p,q)$ be the
coordinates of the most southeast such occurrence.  If at least $1 +
\rr pq\pi$ elements of~$\delta$ are northwest of $(p,q)$, we are done.
On the other hand, if this does not occur, at least $\rr ij\pi - \rr
pq\pi + 1$ elements of~$\delta$ are strictly to the southeast of
$(p,q)$.  By our choice of $(p,q)$, this part of the diagonal
misses~$D(\pi)$.  Hence we obtain a contradiction in view of
Lemma~\ref{lemma:main_obs_unpois} (use the same argument as the one
ending the proof of Proposition~\ref{prop:poison}).%
\end{proof}

\begin{proofof}{Theorem~\ref{thm:neg}}
Assume that $\pi$ is not vexillary.  As the ideals $\IN(A_\pi)$
and~$\IN(B_\pi)$ generated by the diagonal terms of all minors
in~$A_\pi$ and~$B_\pi$ coincide by Proposition~\ref{prop:divisible},
it suffices to prove that $\IN(A_\pi) \neq \IN(I_\pi)$.  The zero set
of $\IN(I_\pi)$ has the same dimension as the zero set of~$I_\pi$,
% because $I_\pi$ and~$\IN(I_\pi)$ have equal Hilbert series,
so it is enough to show that the dimension of the zero set
of~$\IN(A_\pi)$ exceeds that of~$I_\pi$.  Equivalently, it suffices to
prove that $\IN(A_\pi)$ has a component whose codimension is strictly
less than the length of~$\pi$.

The variables corresponding to the $\textcross$ tiles in any pipe
dream~$P$ that poisons~$A_\pi$ generate an ideal $J_P$ that contains
$\IN(A_\pi)$ by definition.  The codimension of~$J_P$ equals the
number of $\textcross$ tiles in~$P$.  The theorem now follows from
Proposition~\ref{prop:poison},
% which says that $C(\pi)$, which has $\ell(\pi)$ crosses,
% poisons~$A_\pi$ but is not a minimal poisoning.
because the number of~$\textcross$ tiles in~$C(\pi)$ is the length
of~$\pi$.%
\end{proofof}

\section*{Acknowledgments}
We would like to thank David Eisenbud, Vic Reiner, Bernd Sturmfels and
Alex Woo for helpful conversations. We thank Frank Sottile for
directing us to the reference \cite{Hodge41}.  AK was supported by an
NSF grant and EM was supported by NSF grant DMS-0304789. This work was
partially completed while AY was an NSERC supported visitor at the
Fields Institute in Toronto, during the 2005 program ``The Geometry of
String Theory''.

%%%%%%%%%%%%%%%%%%%%%%%%%%%%%%%%%%%%%%%%%%%%%%%%%%%%%%%%%%%%%%%%%%%%%%
\end{document}